\documentclass[notitlepage, 11pt]{article}

\usepackage{color}
\usepackage{amsmath, dsfont}
\usepackage{latexsym,amssymb}
\newtheorem{lemma}{Lemma}
\newtheorem{theorem}{Theorem}
\newtheorem{corollary}{Corollary}
\newtheorem{definition}{Definition}
\newtheorem{example}{Example}

\newtheorem{assumption}{Assumption}

\newtheorem{remark}{Remark}

\newcommand{\noi}{\noindent}

\newcommand{\Ra}{\Rightarrow}

\newcommand{\R}{\mathbb{R}}

\newcommand{\p}{\mathbb{P}}
\newcommand{\N}{\mathbb{N}}

\newcommand{\Ir}{\mathbb{I}}

\newcommand{\la}{\lambda}
\newcommand{\sig}{\sigma}

\newcommand{\eps}{\varepsilon}

\newcommand{\kap}{\kappa}

\newcommand{\del}{\delta}

\newcommand{\Del}{\mathnormal{\Delta}}

\newcommand{\Om}{\mathnormal{\Omega}}

\newcommand{\calB}{{\cal B}}

\newcommand{\calF}{{\cal F}}
\newcommand{\calG}{{\cal G}}
\newcommand{\calH}{{\cal H}}

\newcommand{\calM}{{\cal M}}

\newcommand{\calQ}{{\cal Q}}

\newcommand{\calS}{{\cal S}}

\newcommand{\calV}{{\cal V}}

\newcommand{\lan}{\langle}
\newcommand{\lann}{\Bigl\langle}
\newcommand{\ran}{\rangle}
\newcommand{\rann}{\Bigr\rangle}
\newcommand{\lfl}{\lfloor}
\newcommand{\rfl}{\rfloor}

\newcommand{\oo}{\overline}
\newcommand{\skp}{\vspace{\baselineskip}}

\newcommand{\iy}{\infty}

\newcommand{\Gs}{\mathds{G}}

\usepackage{latexsym}
\usepackage{amsmath, amssymb, amsfonts}

\numberwithin{equation}{section}

\title{Fluid limits of G/G/1+G queues under the non-preemptive
earliest-deadline-first discipline\thanks{Research
supported in part by the ISF (Grants 1315/12, 81/09), the US-Israel BSF (Grant 2008466),
and the Technion fund for promotion of research}}

\author{Rami Atar\ \hspace{4em} Anup Biswas  \\[2mm]
	Department of Electrical Engineering\\
Technion--Israel Institute of Technology\\
Haifa 32000, Israel\\
	\and
	Haya Kaspi \\[2mm]
Department of Industrial Engineering and Management\\
Technion--Israel Institute of Technology\\
Haifa 32000, Israel
	}

\date{January 29, 2014}

\begin{document}

\maketitle

\begin{abstract}
A single-server queueing model is considered with customers that have deadlines.
If a customer's deadline elapses before service is offered, the customer
abandons the system (customers do not abandon while being served).
When the server becomes available, it offers service to the customer having earliest deadline
among those that are in the queue.
We obtain a fluid limit of the queue length and abandonment processes and for the occupation
measure of deadlines, in the form of measure-valued processes.
We characterize the limit by means of a Skorohod problem in a time-varying domain,
which has an explicit solution.
The fluid limits also describe a certain process called the {\it frontier},
that is well known to play a key role in systems operating under this scheduling policy.

\skp

\noi{\bf AMS subject classifications:}\, 60K25, 60G57, 68M20

\skp

\noi{\bf Keywords:}\, due dates, G/G/1+G, earliest deadline first,
fluid limits, law of large numbers, measure-valued processes
\end{abstract}

\section{Introduction}

In recent years there has been growing interest in queueing models with reneging customers.
An important family of such models is one where arriving customers
have deadlines that are known to the decision maker,
and they renege the system if service does not start (or end) by the time of the deadline.
Well known examples of such systems are time-sensitive computer networks,
real-time control systems and telecommunication systems for transferring audio or video signals.
In case of transmission of voice or video over a packet-switched network,
the packets have to reach their destination within a certain time of their transmission, and
otherwise they are useless to the receiver and considered lost.
The performance of such a system is measured by its ability to meet the deadline constraints.
Therefore it is desired to schedule the service so as to maximize the fraction of customers
served before their deadline expires.

In this article we consider a single-server queueing model with customer deadlines,
operating under the earliest-deadline-first (EDF) scheduling discipline,
according to which, when the server completes a job it accepts the customer
whose deadline is the earliest among all customers that are still in the system.
In addition to the {\it deadline} we will sometimes talk of a customer's {\it lead time}
defined as the difference
\begin{center}
lead time of a customer $=$ deadline of the customer $-$ current time.
\end{center}
According to this terminology, a customer reneges if its lead time becomes negative
before it starts being served. Moreover, service is provided to customers with positive
lead time, who are prioritized according to their lead times.

There are several variants of the EDF policy. The deadlines could refer to the time
of beginning of service, in which case customers do not renege while in service.
We refer to this policy as EDF-b.
Or it could refer to the end of service, i.e., customers will renege
whenever their deadlines elapses irrespective of their status, a variant which we call
EDF-e. One can also consider preemptive and non-preemptive policies.
`Preemptive' refers to switching from serving a customer to serving
a newly arrived customer if it has a deadline that is earlier
than the one in service; `non-preemptive' means that service can not be interrupted.

Various queueing models employing EDF policies have been studied by exact analysis
as well as under scaling limits. Results on exact analysis include the following.
Panwar and Towsley \cite{panwar-towsley}
prove the optimality of the preemptive EDF-e policy within the
class of work-conserving policies for the G/M/1+G queue.
Kruk et al.\ \cite{kruk-lehoc-ram-shre} study a G/G/1+G queue
and use the amount of {\it reneged work} as a performance measure,
showing that the preemptive EDF-e policy is optimal.
Panwar, Towsley and Wolf \cite{panwar-towsley-wolf}
show that the non-preemptive EDF-b policy maximizes
the expected number of customers to meet their deadlines,
within the class of work-conserving, non-preemptive policies for the M/G/1+G system.
Moyal \cite{Moyal}
studies stability and optimality properties of non-preemptive EDF-b,
under various cost functions, for single-server queueing systems.

The first to consider the EDF policy for the G/G/1+G queue under scaling limits is the work by
Doytchinov, Lehoczky and Shreve \cite{Doy-Leh-Shre}, which studies the preemptive EDF-e
policy under heavy traffic diffusion approximations.
It is assumed in this work that the system serves all the customers, including those who missed
their deadlines. The approach is to work with measure-valued processes:
the {\it queue-length measure-valued process} ({\it queue measure} for short),
corresponding to
the empirical measure of lead times of all customers in the system, as it varies in time,
and the {\it workload measure-value process} ({\it workload measure} for short) corresponding to the empirical measure of
workload carried by each of the customers in the system.
It is established that these measure-valued processes converge, and their limits
are identified. This work was extended
by Kruk et al.\ in \cite{kruk-lehoc-ram-shre} to treat a system where customers
who miss their deadline do renege,
for patience with strictly positive lower bound. In this work it is proved that
the workload, queue-length and reneged work processes converge, and their
limit is identified.
In both \cite{Doy-Leh-Shre} and \cite{kruk-lehoc-ram-shre},
the {\it frontier} process plays a key role in identifying the limit of the above processes.
The frontier at time $t$ is defined there as the largest
lead time among all customers who have ever been in service at or prior to $t$.
The work which is closest to our present paper is
Decreusefond and Moyal \cite{Dec-Moyal},
which considers an M/M/1+G queue with the non-preemptive EDF-b policy
and establishes the fluid limit
of the queue measure and {\it counting process of reneged customers}
({\it reneging process} for short).
The main idea of \cite{Dec-Moyal} is to use the Markov property in the evolution of
the point measure-valued process, assigning point mass at the lead times of the customers
either in queue or discarded. A key assumption is
the convergence of the least lead time process to a deterministic continuous path, as assumption
that is verified in the same paper
only in the case of M/M/1+D, where ``+D'' stands for deterministic lead times.

In the present article we aim at fluid limits for the more general setting of the G/G/1+G queue,
operating under the non-preemptive EDF-b policy. We obtain fluid limits of the queue length
and reneging processes, as well as the measure-valued queue-length process.
Our work does not extend, strictly speaking, the work \cite{Dec-Moyal}, in the sense that
our assumptions exclude the case of deterministic lead times.
Our results also address the frontier process, although
we adopt a slightly different definition of the frontier, which we take to
be, at time $t$, the largest
lead time among all customers who have been at the head of the queue at or prior to $t$.
We prove that the frontier converges to a deterministic continuous function,
and that the queue measure is governed by the frontier, in the critical and super-critical
cases.

At the crux of our proof is an argument relating the queue-length, the reneging and
the {\it potential queue-length} processes.
The latter quantity, introduced in Kang and Ramanan \cite{kang-ramanan}, is defined as
the queue-length in a model that has no server, but has
arrival and reneging as in the original model. In other words,
it is the queue-length that would be obtained in our model if services are switched off.
The three processes alluded to above
satisfy a relation that, in the limit, is characterized by means of
a one-dimensional Skorohod Problem with a time-varying barrier.
In this relation, the limiting potential queue-length process
serves as the barrier, while the limiting queue-length and reneging processes act as the solution.
The potential queue-length can be expressed directly in terms of the primitive processes,
and its law-of-large-numbers limit is found in a straightforward way.
As a result, a characterization,
and in fact an explicit formula for the queue-length and reneging limits follow immediately.
This relation is closely related to an idea from
\cite{kruk-lehoc-ram-shre}, where it is shown, under diffusion scaling and
a preemptive EDF-e policy, that the queue-length (as well as the workload)
is governed by a {\it doubly}-reflected Brownian motion, on a {\it fixed} interval.
The proof in \cite{kruk-lehoc-ram-shre} relies on the analysis of a measure-valued process
referred to as the \textit{reference workload}, and applies an argument
that the actual workload process asymptotically merges with the reference workload process.
In the present paper our approach of obtaining the Skorohod problem relation is much more direct,
and in particular does not rely on the underlying measure-valued process.
Rather, the relation is
a direct consequence of an analogous relation satisfied by the three prelimit processes.
This way we can obtain the limits of the three one-dimensional processes without relying on
limit results for the measure-valued processes.
The analysis of the measure-valued processes is carried out only afterwards,
and builds on this result.

Measure-valued processes have been useful in modeling and analyzing
queuing systems beyond the EDF policy in various works, including the following.
Gromoll, Puha and Williams \cite{grom-puha-will} and
Zhang, Dai and Zwart \cite{zhang-dai-zwart} use measure-valued processes
to obtain fluid limits for processor sharing models and limited
processor sharing models, respectively.
Down, Gromoll and Puha \cite{down-grom-puha}
consider a single-server queue with the shortest-processing-time-first policy, establishing the
fluid limit of the workload measure.
The diffusion limit for a critically loaded
queue with general service and abandonment time distributions, under
FIFO discipline, is studied by Glynn and Ward \cite{glynn-ward}
and Reed and Ward \cite{reed-ward}.
Kaspi and Ramanan \cite{kaspi-ramanan} use measure-valued processes
to model the G/G/N queue with $N\to\iy$, obtaining their fluid limits.
The paper by Kang and Ramanan \cite{kang-ramanan}, mentioned above,
obtains fluid limit for the G/G/N+G queue under a similar scaling.
Atar, Kaspi and Shimkin \cite{atar-kas-shim} use the methods from \cite{kaspi-ramanan} and
\cite{kang-ramanan} to establish the fluid limit of the multi-class G/G/N+G
under a priority policy. Biswas \cite{biswas} establishes the fluid limit a many-server
queueing system with state-dependent service rates.

\skp

\noi{\bf Notation:}
The following notation will be used throughout this paper.
For $x\in\R$, $\lfl x\rfl$ denotes the largest integer less than or equal to $x$.
For $x,y\in\R$, the maximum (minimum) is denoted by $x\vee y$ (resp., $x\wedge y$), and
$x^+=x\vee 0$.
The symbol $\del_x$ denotes the point mass at $x\in\R$.
We denote $\R_+=[0, \infty)$.
For any $A\subset [0,\infty)$, we define $A^\eps=\{x\geq 0: \inf_{a\in A}|x-a|<\eps\}$. For any $x\in\R$, the sets $(x, \infty), [x, \infty)$
will be denoted by $C_x, \bar C_x$ respectively. Characteristic function of a set $A$ is denoted by $\Ir_A$.
For a topological space $\calS$,
$C_b(\calS)$ denotes the set of real-valued bounded, continuous maps on $\calS$, while $\mathcal{B}(\calS)$ denotes the corresponding Borel $\sigma$-field.

The space of non-negative finite Borel measures on $[0,\infty)$ is denoted by $\calM$. For any $\mu\in\calM$ and Borel measurable function $g$ on $[0,\infty)$, denote
$\langle g, \mu\rangle=\int g d\mu$.
Endow $\calM$ with the Prohorov metric given by
$$\rho(\mu_1,\mu_2)=\inf\{\eps>0 : \mu_1(A)\leq \mu_2(A^\eps)+\eps, \mu_2(A)\leq \mu_1(A^\eps)+\eps,
\mbox{for all closed}\ A\subset[0,\infty)\}.$$
It is well known that $(\calM, \rho)$ is a Polish space \cite[Appendix]{delay-verejones}. Also this topology is equivalent
to the weak topology on $\calM$, characterized by $\mu_n\rightarrow \mu$
if and only if
\begin{center}
$\langle f, \mu_n\rangle\to\langle f, \mu\rangle$ for all $f\in C_b(\R_+)$.
\end{center}
For any function $f:\R_+\to\R$ and $\del>0$, $osc_\del(f)$ denotes the $\del$-oscillation of $f$ given by
$$osc_\del(f)=\sup\{|f(s)-f(t)|\ :\ |s-t|\leq\del\}. $$
By $\Delta f(t)$ we denote the jump of $f$ at $t$ i.e., $\Delta f(t)=f(t)-f(t-)$.
For a Polish space $\calS$, denote by $D_{\calS}(\R_+)$
and $C_\calS(\R_+)$ the spaces of functions $\R_+\to\calS$ that are
right continuous with finite left limits, and respectively,
continuous. Endow the space $D_{\calS}(\R_+)$
with the Skorohod $J_1$ topology. All the random variables introduced in this paper
are assumed to be defined on a common probability space,
$(\Om, \calF, \p)$. We use the symbol `$\Ra$' to denote the convergence
in distribution of a sequence of random variables.
For a cumulative (sub-) distribution function $G$ we write $\oo G$ for
$G(\iy)-G(x)$ (in the case where $G$ corresponds to a probability, this is $1-G$).

The organization of the paper is as follows. The next section introduces the model and
states
the main results. Section 3 presents the proofs of the main theorems along
with subsequent lemmas.

\section{The queueing model and main results}

We begin by giving a precise description of the model.
We consider a sequence of systems, indexed by $N\in\N$,
having a single server and a buffer of infinite room for customers waiting to be served.
Each customer has a single service requirement and leaves the system when his job is completed.
The arrival process for the $N$-th system,
denoted by $E^N$, is assumed to be a renewal process with mean
inter-arrival time $\frac{1}{\la^N}$, where $\la^N>0$ are given parameters.
These processes are obtained by time acceleration of
a single processes, namely $E^N(t)=E^0(\la^Nt)$, $N\ge1$, where $E^0$ has inter-renewal mean 1.
Denote by $a^N_i$ the arrival time of the $i$-th customer in the $N$-th system, $i=1,2,\ldots$
The service requirements are assumed to be generally distributed with
{\it mean service time} $\frac{1}{\mu^N}$, where $\mu^N>0$ are given parameters.
To the customers that are in the system at time $0$ we associate their arrival times $a^N_i, i=-X^N(0)+1, \ldots, 0$, satisfying $a^N_{-X^N(0)+1}\leq a^N_{-X^N(0)+2}\leq\cdots\leq a^N_0$, where $X^N(0)$ denotes the number of customers in the system at time zero.
An incoming customer has an {\it initial lead time}, expressing how long after arrival he
will wait in the queue before abandoning the system.
For each $N$, the initial lead times are given by
a positive i.i.d.\ random variables $\{u_i^N\}$. For the $i$-th customer, there is an
associated stochastic process called the {\it lead time} $L_i^N(t)$ given by
the initial lead time minus the time spent in the system. That is,
\begin{equation}\label{r01}
L_i^N(t)=u_i^N+a^N_i-t,\qquad t\ge0.
\end{equation}
We assume that the sequence of initial lead times,
the arrival process and the service requirements are mutually independent.
After arrival, a customer could be in one of three states: in the queue, in service, or
outside the system.
Upon arrival, customers wait in the queue
if the server is busy, or otherwise immediately start service.
A customer waiting in the queue
{\it reneges} the system once his lead time becomes less than or
equal to zero. Thus, at any time, all customers present in the queue have positive lead times.
Customers are served according to the non-preemptive EDF-b policy, described as follows.
When the server completes a job and the queue in not empty,
it starts providing service to the customer having the smallest lead time
(equivalently, smallest positive lead time) among the customers in the queue.
Note that by this description the policy
is non-preemptive and non-idling. By $X^N(t)$ and $Q^N(t)$, we denote the number of customers
in the system, and respectively, in the queue at time $t$.
By $R^N(t)$ we denote the number of customers having reneged in $[0, t]$.
We assume $R^N(0+)=0$. $R^N$ is called the reneging process.

Let $S^N(t)$ be the number of jobs completed by the time the server has been busy for
$t$ units of time. By assumption it is
a renewal process with mean inter-renewal time $\frac{1}{\mu^N}$, and we assume that
$S^N(t)=S^0(\mu^N t)$, $N\ge1$, where $S^0$ has inter-renewal mean 1.
The number of service completions of jobs by time $t$ is thus given by
\begin{equation}\label{101}
D^N(t)=S^N\Big(\int_0^t\Ir_{\{X^N(s)>0\}}ds\Big).
\end{equation}
Hence we have
\begin{align}
X^N(t)&= X^N(0)+E^N(t)-D^N(t)-R^N(t),\label{102}
\\
Q^N(t)&=(X^N(t)-1)^+.\label{103}
\end{align}

Denote by $\tau^N_i,  i\geq -X^N(0)+1$, the time when the $i$-th customer starts its
service. We define $\tau^N_i$ to be infinity if the $i$-th customer reneges the system.
Now we define the {\it queue-lead-time measure-valued process},
or the {\it queue measure} for short, by

\begin{equation}\label{queue-mes}
\calQ^N_t(B)=\sum_{ i\le E^N(t), i\neq 0:\ t<\tau^N_i<\iy}\Ir_B(L_i^N(t)),\qquad t\ge0,
\end{equation}
where $L_i^N(t)=L_i^N(0)-t$ for $i\le0$ and $B\in \calB((0, \infty))$ .
Then for $B\subset(0,\iy)$ , $\calQ^N_t(B)$ represents the number of customers in the queue at time $t$
whose lead times, at time $t$, are in $B$.
We define $F^N(0)$ as the lead time of the customer at the head of the queue at time $0$
or $0$ if the queue is empty.
Define
\begin{equation}\label{lead-Z}
Z^N(t)=
\begin{cases}
\sup\{x\geq 0\ : \calQ^N_t([0, x))=0\} & \mbox{if}\ Q^N(t)\neq 0,
\\
0 & \mbox{otherwise.}
\end{cases}
\end{equation}
For $Q^N(t)\neq 0$, $Z^N(t)$ is the lead time of the customer at the head of
the line.
We also define the {\it frontier}
\[F^N(t)=\Big(\max_{s\leq t; Q^N(s)\neq 0}\{Z^N(s)-t+s\}\Big)\vee (F^N(0)-t).
\]
By this definition, $F^N(t)$ equals
the largest lead time, at time $t$, of any customer who has ever been at the head of
the queue, whether still present in the system or not, or $F^N(0)-t$ if this is
larger than the former.
The {\it current lead time} $C^N(t)$ is defined by
\[C^N(t)=\left\{\begin{array}{ll}
Z^N(t) & \mbox{if} \ Q^N(t)\neq 0,
\\
F^N(t) & \mbox{otherwise.}
\end{array}
\right.
\]
Then  $C^N(t)$ equals the lead time, at time $t$,
of the customer who is in the head of the queue, or $F^N(t)$ if the queue is empty at time
$t$.

By the above description we see that $\calQ^N_t$ is a measure on $(0, \iy)$ with
$\calQ^N_t[0, C^N(t))=0$ for all $t$.
Moreover, it is easy to see that $F^N$ and $C^N$ have RCLL sample paths.
The process $F^N$ may take negative values, but $F^N(t)>0$ if
$Q^N(t)>0$. Finally, $F^N(0)=C^N(0)\geq 0$.

The idea of analyzing the frontier process in order to obtain scaling limits has been
employed very successfully for systems that operate under priority based
on deadlines, as in \cite{Dec-Moyal}, \cite{Doy-Leh-Shre} and \cite{kruk-lehoc-ram-shre},
as well as on job size, as in \cite{down-grom-puha}. The papers
\cite{Doy-Leh-Shre} and \cite{kruk-lehoc-ram-shre} use frontier process to obtain diffusion limits and
the same is used in \cite{Dec-Moyal} to obtain fluid limits when the customers have deterministic deadlines.
Specifically,
\cite{Doy-Leh-Shre}, \cite{kruk-lehoc-ram-shre} consider queuing systems under preemptive EDF and the frontier is defined
as the largest lead time among the customers who have ever been in service. Under diffusion setting it is shown that the re-scaled customer population with
lead times in $[C^N, F^N]$ disappears asymptotically. Here we modify the definition of frontier so that it is suited to a non-preemptive policy. We also show that $\frac{1}{N}\calQ^N([C^N, F^N])$ tends to $0$ (Lemma \ref{lem-5}).
Under the first order scaling, the limit of the frontier processes is identified in \cite{Dec-Moyal} when the deadlines are deterministic, in which case
one has $F^N=C^N$.
Let us mention that our assumptions do not allow the deadlines to be deterministic.
\skp

Let $\nu^N$ be the distribution of $u^N_1$ and $G^N(\cdot)$ its cumulative distribution. Define
$y^N_{\max}=\inf\{x:\ G^N(x)=1\}$. We shall consider a sequence of such queues indexed by $N$, and make the following assumptions.

\begin{assumption}\label{cond1}
As $N\to\iy$
\begin{itemize}
\item[(i)] $\frac{\la^N}{N}\to\la>0$,
\item[(ii)] $\frac{\mu^N}{N}\to\mu>0,$
\item [(iii)] $u^N_1$ converges in distribution to a random variable $u$ with distribution $\nu$ and cumulative distribution function $G$.
We assume that $G$ is strictly increasing in $(0, y_{\max})$ where $y_{\max}=\inf\{x:\ G(x)=1\}$ and $G$ is continuous at $0$ with $G(0)=0$.
\end{itemize}
\end{assumption}

Let
\begin{equation}\label{y-cond}
 y^*=\sup\{y< y_{\max}\ :\ \la G(y)<\mu\}.
\end{equation}

\begin{assumption}\label{G-cond}
 For $\la>\mu$, we assume that $G$ is continuous at $y^*$.
\end{assumption}

We assume the following condition on the initial queue measure
\begin{assumption}\label{initial-meas}
There exists a pair $(\calQ_0, F(0))\in\calM\times [0, y^*]$ so that
\begin{itemize}
\item[(i)]
$\calQ_0([0, F(0)])=0$ and $G_0(x)=\calQ_0([0, x])$ is a continuous cumulative distribution function.
\item[(ii)] For all $x>F(0)$, $\calQ_0([0,x])>0$.
\item [(iii)] $(\frac{1}{N}\calQ^N_0, F^N(0))\Ra(\calQ_0, F(0))$  as $N\to\iy$.
\end{itemize}
\end{assumption}

\begin{remark}\label{rem1}
For $\la\leq\mu$, one has
$0\leq\frac{1}{N}Q^N\leq\frac{1}{N}(Q^N(0)+E^N(t)-D^N(t))$.
For large $N$ this is close to $\frac{1}{N}(Q^N(0)+N\la t-N\mu t)$
which is smaller than $\frac{1}{N}Q^N(0)$, provided that $[0, t]$ is contained in one busy period. Thus to obtain a non-zero fluid limit one needs to assume that $\frac{1}{N}Q^N(0)$ converges to $Q(0)>0$.
\end{remark}
Recall our notation $\oo G=1-G$. Let $\oo G_0(\cdot)=\calQ_0(\cdot, \iy),$ and
\begin{equation}\label{301}
H(x, t)=\oo{G}_0(x+t)+\la\int_0^t\oo{G}(x+t-s)ds.
\end{equation}
It is not hard to see that $H$ satisfies the following transport equation
\begin{equation*}
\begin{array}{ll}
\partial_t H(x, t)&=\partial_x H(x, t)+\lambda \oo{G}(x),\quad \ \text{in}\ \R_+\times(0, \iy),
\\
H(x, 0) &= \oo{G}_0(x), \quad \ \text{in}\ \R_+.
\end{array}
\end{equation*}
An approach that uses the transport equation in this context
was introduced in \cite{Dec-Moyal}, where it plays a central role in the analysis.
Our analysis does not use the equation but works directly with the definition \eqref{301}.

By Assumptions \ref{cond1} and \ref{initial-meas} for $T>0$, $H: \R_+\times[0, T]\to\R_+,$ is continuous and for a fixed $t\in(0, T]$, $H(\cdot, t)$
is a continuous,
 decreasing function from $\R_+$ onto its range. In fact, $H(\cdot, \cdot)$ is uniformly continuous on
$\R_+\times[0, T]$. To see this it is enough to show that $\int_0^\cdot \oo{G}(\cdot+s)ds$ is uniformly continuous. Let $0\leq t_1<t_2\leq T$ and $x_1, x_2\in\R_+$. Then
\begin{align*}
\Big|\int_0^{t_2} \oo{G}(x_2+s)ds- \int_0^{t_1} \oo{G}(x_1+s)ds\Big| &=
\Big|\int_{x_2}^{x_2+t_2} \oo{G}(s)ds- \int_{x_1}^{x_1+t_1} \oo{G}(s)ds\Big|
\\
&\leq 2|x_2-x_1|+(t_2-t_1).
\end{align*}

Now we introduce a {\it Skorohod problem} on a time-varying domain.
\begin{definition}\label{def1}
Let $h\in C_\R(\R_+)$ be fixed. Given $\psi\in C_\R(\R_+)$, a pair $(\phi,\eta)\in C_\R(\R_+)^2$
is said to solve the Skorohod Problem (SP) for $\psi$ on the time-varying domain
$(-\iy,h(\cdot)]$, if
\begin{enumerate}
\item $\phi(t)=\psi(t)-\eta(t), \quad t\geq 0$,
\item $\phi(t)\leq h(t),\quad t\geq0$,
\item $\eta$ is non-negative, non-decreasing and
$\int_0^\cdot \Ir_{\{\phi(s)<h(s)\}}d\eta(s)=0$.
\end{enumerate}
\end{definition}
In fact, there is always a unique pair $(\phi, \eta)$ solving the SP for a given $\psi$,
given by (\cite{burdzy-kang-ramanan})
\[
\eta(t)=\sup_{s\in[0, t]}(\psi(s)-h(s))^+,
\qquad
\phi(t)=\psi(t)-\eta(t),\qquad t\ge0.
\]
Henceforth we denote the {\it Skorohod map} by $\Gamma_h$ i.e, $\Gamma_h(\psi)=\phi$.
We will be interested in the case where $h=H(0,\cdot)$ and
$\psi(t)=Q(0)+(\la-\mu)t$. The solution to the corresponding
SP in this case is given by
\begin{equation}\label{105}
\eta(t)=\sup_{s\in[0, t]}(\psi(s)-H(0,s))^+,
\qquad
\phi(t)=\psi(t)-\eta(t),\qquad t\ge0.
\end{equation}

\skp

We now state our main results
\begin{theorem}\label{main2-th}
Under Assumptions \ref{cond1}--\ref{initial-meas} and $\la\geq\mu$,
let $\psi(t)=Q(0)+(\la-\mu)t$ and let $(\phi,\eta)$ be
given by \eqref{105}.
Then $(\frac{1}{N}Q^N,\frac1NR^N)\Ra(\phi,\eta)$ as $N\to\iy$.
\end{theorem}

Define $\chi:\R_+\times[0, \iy)\to\R_+$
as follows
\begin{equation}\label{inv}
\chi(x, t)=\inf\{y\geq 0 :\ H(y,t)\leq x \}.
\end{equation}
Define $y_{0\min}=\inf\{y\ :\ G_0(y)>0\}$ and
$y_{0\max}=\sup\{y\ :\ G_0(y)<1\}$. Then $[y_{0\min}, y_{0\max}]$
is the support of $\calQ_0$.

\begin{theorem}\label{main-th}
Suppose that $\la\geq\mu$ and
Assumptions \ref{cond1}--\ref{initial-meas} are satisfied. Assume further that $\oo{G}_0$
is strictly decreasing in
$[y_{0\min}, y_{0\max}]$. If $y_{0\min}=0$ and  $y_{0\max}\geq y_{\max}$(recall that $y_{\max}=\inf\{x:\ G(x)=1\}$), then
\begin{itemize}
\item[(a)] $F^N\Ra F$ as $N\to\iy$, where
\begin{equation}\label{106}
F(t)=\chi(\phi(t), t), \quad t\geq 0,
\end{equation}
and $\phi$ is as in Theorem \ref{main2-th}.

\item[(b)] Let $\calQ$ be the $\calM$-valued RCLL path defined by
$$\calQ_t(B)= \calQ_0(B\cap[F(t),\iy)+t)+\la\int_0^t \nu(B\cap[F(t),\iy)+t-s) ds,$$
where $\nu$ is defined in Assumption \ref{cond1}, $B\in\calB([0,\iy))$ and $A+t=\{x+t:\ x\in A\}$. Then
$$ \frac{1}{N}\calQ^N\Ra\calQ\quad \text{as}\quad N\to\iy.$$
\end{itemize}
For general $y_{0\min}$ and $y_{0\max}$ the above convergence holds on compact subsets of $(0, T)$ for $T>0$.
\end{theorem}

\begin{theorem}\label{main3-th}
Suppose that Assumption \ref{cond1} and \ref{initial-meas} hold and that $\la<\mu$.
Let  $\bar T=\inf\{t\geq 0\ : \phi(t)=0\}$ and define
\[\bar\phi(t)=\left\{
\begin{array}{ll}
\phi(t) & \text{if}\ t\leq\bar T,
\\
0 & \mathrm{otherwise}.
\end{array}
\right.
\]
Then $\frac{1}{N}Q^N\Ra\bar\phi$ as $N\to\iy$.
\end{theorem}

\begin{remark}
The definition of the frontier used in this article is expected to be useful in other scenarios as well. For example, one can consider
a queueing model with many servers operating under the EDF policy, in which limits
are taken in a fashion similar to \cite{kaspi-ramanan} where the number of servers
tends to infinity.
Similar results are expected there
for an overloaded system when the job requirement is exponentially distributed.
\end{remark}

\begin{example}
We demonstrate that the ingredients of the formula for the frontier
process limit \eqref{106} can sometimes be computed explicitly.
Consider a system with $\la=1$, $\mu\le 1$ and $\theta>0$.
Assume $\bar G_0(x)=e^{-x}$.
Then $G(x)=1-e^{-\theta x}$, $Q(0)=1$ and
$$
H(x, t)=e^{-x-t}+\frac{1}{\theta}e^{-\theta x}(1-e^{-\theta t}).
$$
Therefore $\psi(t)=1+(1-\mu)t$. To compute $(\phi, \eta)$ from \eqref{105}
we distinguish three cases.

\noi{Case 1:} $\mu\leq 1$ and $\theta\geq 1$. Then
\begin{align*}
\phi(t)&= \frac{1}{\theta}+e^{-t}-\frac{1}{\theta}e^{-\theta t},
\\
\eta(t) & = 1-\frac{1}{\theta}+(1-\mu)t-e^{-t}+\frac{1}{\theta}e^{-\theta t}.
\end{align*}

\noi{Case 2:} $\mu=1$ and $\theta<1$. Then $\phi(t)=1$ and $\eta(t)=0$ for all $t$.

\noi{Case 3:} $\mu<1$ and $\theta<1$. If $1-\mu\geq (1-\theta)\theta^{\frac{\theta}{1-\theta}}$ then
$(\phi, \eta)$ is the same as in Case 1. If $1-\mu< (1-\theta)\theta^{\frac{\theta}{1-\theta}}$
then there are constants $a_1, a_2$ (in fact, $a_1$ is the smallest zero of $f(s)=1-\mu+e^{-s}-e^{-\theta s}$ and $a_2$ is the largest
point satisfying $\int_0^{a_1}f(s)ds=\int_0^{a_2}f(s)ds$) such that
\[\eta(t)=\left\{
\begin{array}{lll}
1-\frac{1}{\theta}+(1-\mu)t-e^{-t}+\frac{1}{\theta}e^{-\theta t} & \text{for}\quad t\leq a_1,
\\
\eta(a_1) & \text{for}\quad t\in[a_1, a_2],
\\
1-\frac{1}{\theta}+(1-\mu)t-e^{-t}+\frac{1}{\theta}e^{-\theta t} & \text{for}\quad t\geq a_2,
\end{array}
\right.\]
and $\phi(t)=\psi(t)-\eta(t)$.
\end{example}

\section{Proofs}
This section is devoted to the proof of
Theorem \ref{main2-th}, \ref{main-th} and \ref{main3-th}. First we prove some lemmas that will be
used to prove the main results.

Define $\bar E^N=\frac{1}{N}E^N$. It follows from
Assumption \ref{cond1}(i), $\bar E^N(t)\to\la t$ in probability as $N\to\iy$,
u.o.c.
Hence, given $T$, there exists a sequence $\eps_A(N)\to 0$ such that
\begin{equation}\label{201}
\p(\sup_{t\in[0,T]}|\bar E^N(t)-\la t|\leq\eps_A(N))\geq 1-\eps_A(N).
\end{equation}
Similarly, defining $\bar S^N(t)=\frac{1}{N}S^N(t)$, there exists a sequence
$\eps_S(N)\to 0$ such that
\begin{equation}\label{212}
\p(\sup_{t\in[0,T]}|\bar S^N(t)-\mu t|\leq\eps_S(N))\geq 1-\eps_S(N).
\end{equation}
Let $\Om_A(N)$ [resp., $\Om_S(N)$] denote the event indicated on the l.h.s.\ of \eqref{201}
[resp., \eqref{212}].
Recall that $C_x=(x, \iy)$ and $\bar C_x=[x, \iy)$ for $x\in\R_+$. We prove the following lemma for general cumulative distribution function $G$. Recall from \eqref{r01} the notation $L^N_i(t)$
for the lead time of the $i$-th customer at time $t$.
\begin{lemma}\label{lem-2}
Let $\eps, \eta>0$ be given. Then for $T>0$,
$$
\limsup_{N\to\iy}\p\Big(
\sup_{a\in\R_+}\sup_{f\in\{\Ir_{C_a}, \Ir_{\bar C_a}\}}\sup_{t\in[0, T]}
\Big|\lan f, \calG^N_t
\ran-\la\int_0^t\oo{G}(a+t-s)ds\Big|\geq\eps\Big)\leq \eta,
$$
where
$$\calG^N_t= \frac{1}{N}\sum_{i=0}^{E^N(t)}\del_{L^N_i(t)}.$$
\end{lemma}

\noi{\bf Proof:} Write $G$ as
$G=\mathds{G}_c+\Gs_d$,
where
\[
\Gs_d(x)=\sum_{y\leq x}\Del G(y),\quad \Del G(y)=G(y)-G(y-),
\]
and $\Gs_c$ is a continuous function. It is also clear that $G$ being non-decreasing, $\Gs_c$ is also non-decreasing. For $x, y\in\R_+, x\geq y$,
$$ \Gs_c(x)-\Gs_c(y)=[G(x)-G(y)]-[\sum_{y<z\leq x}\Delta G(z)]\geq 0.$$
Recall our notation
$\oo{\Gs}_c(x)=\Gs_c(\iy)-\Gs_c(x)$ and $\oo{\Gs}_d=\Gs_d(\iy)-\Gs_d$. Then $\oo{G}=\oo{\Gs}_c+\oo{\Gs}_d$
and $\sup_{x\in\R_+}|\Gs_d(x)|\leq 1$.
First we show that for $a\in\R_+$ and $\eps>0$,
\begin{align}
\label{202}
&\lim_{N\to\iy}\p\Big(\sup_{f\in\{\Ir_{C_a}, \Ir_{\bar C_a}\}}\sup_{t\in[0, T]}\Big|\lan f, \calG^N_t\ran\notag\\
&-\int_0^t\oo{\Gs}_c(a+t-s)d\bar E^N(s) -\la\int_0^t\oo{\Gs}_d(a+t-s)ds\Big|\geq\eps\Big)=0.
\end{align}
Since $G$ is non-decreasing it can have at most countably many jumps in $[0, \iy)$, say $\{g_i\}$.
Now given $\eps>0$ we can find a positive integer $N_\eps$ such that $\sum_{j\geq N_\eps}\Delta G(g_j)\leq\frac{\eps}{16(\la T+1)}$.
We order the set $\{g_1,\ldots, g_{N_\eps}\}$ as $\{\tilde g_1, \ldots, \tilde g_{N_\eps}\}$ so that $\tilde g_1<\ldots<\tilde g_{N_\eps}$.
Now define a function $\Gs_{d, \eps}$ on $[0, \iy)$ as
\[\Gs_{d,\eps}(x)=\left\{\begin{array}{lll}
\oo{\Gs}_d(0)\ & \text{if}\ x\in [0, \tilde g_1),
\\
\oo{\Gs}_d(\tilde g_i) & \text{if}\ x\in [\tilde g_{i}, \tilde g_{i+1}), 1\leq i\leq N_\eps-1,
\\
\oo{\Gs}_d(\tilde g_{N_\eps})& \text{if}\ x\in [\tilde g_{N_\eps}, \iy).
\end{array}
\right.
\]
Then it is easy to see that
\begin{equation}\label{1000}
\sup_{x\in\R_+}|\oo{\Gs}_d(x)-\Gs_{d, \eps}(x)|\leq \frac{\eps}{16(\la T+1)}.
\end{equation}
We choose $\del>0$ and a partition $\{t_i\}_{i=0}^{i=k(t)}$ of $[0, t]$ such that $\max_{1\leq i\leq k(t)}|t_i-t_{i-1}|\leq \del$
and $\sup_{t\in[0, T]}k(t)=k_\del<\iy$. We can choose the partition to satisfy
$2\del<\min_{i\neq j}|\tilde g_i-\tilde g_j|$ and
for all $i$ such that $\tilde g_i\le a+t$,
$\tilde g_i=a+t-t_j$ for some $j\le k(t)$.
Then
\begin{equation}\label{260}
\la\int_0^t\Gs_{d, \eps}(a+t-s)ds=\la\sum_{m=1}^{k(t)}\Gs_{d, \eps}(a+t-t_{m})(t_m-t_{m-1}).
\end{equation}
If we choose $N$ large enough then on $\Om_A(N)$ we have
\begin{equation}\label{211}
\sup_{t\in[0,T]}\bar E^N(t)\leq \la T+1.
\end{equation}
Define
$$ \calV=\{\Ir_{C_x}\ :\ x\in\R_+\}\cup\{\Ir_{\bar C_x}\ :\ x\in\R_+\}.$$
Then defining $M=\la T+1,$
we have, using \cite[Appendix B]{zhang-dai-zwart},
$$\limsup_{N\to\iy}\p\Big(\max_{m\in[0, NM]}\sup_{\ell\in[0, M]}\sup_{f\in\calV}\Big|\Big\lan f, \frac{1}{N}\sum_{i=1+m}^{m+\lfl N\ell\rfl}
\del_{u_i^N}\Big\ran-\ell \lan f,\nu^N \ran\Big|>\eps\Big)<\eta. $$
As a result there exists a sequence $\{\eps_G(N)\},\eps_G(N)\to 0$ as $N\to\iy$,
such that for
\begin{equation}\label{90}
\Om_G(N)=\Big\{\max_{m\in[0, NM]}\sup_{\ell\in[0, L]}\sup_{f\in\calV}\Big|\lan f,
\frac{1}{N}\sum_{i=1+m}^{m+\lfl N\ell\rfl}\del_{u_i}\ran-\ell\lan f,\nu^N\ran\Big|
\leq\eps_G(N)\Big\}
\end{equation}
we have $\p(\Om_G(N))\geq 1-\eps_G(N)$.
Hence for large $N$, on $\Om_A(N)\cap\Om_G(N)$, we have for $f=\Ir_{C_a}$,
\begin{align*}
&\lan f, \calG^N_t\ran - \int_0^t\oo{\Gs}_c(a+t-s)d\bar E^N(s)-\la\int_0^t\oo{\Gs}_d(a+t-s)ds \\
&\leq \frac{1}{N}\sum_{i=0}^{E^N(t)} \lan \Ir_{C_{a+t-a^N_i}},\del_{u_i}\ran
-\int_0^t\oo{\Gs}_c(a+t-s)d\bar E^N(s) -\la\int_0^t\Gs_{d, \eps}(a+t-s)ds+\frac{\eps}{16}
\\
& =\sum_{m=0}^{k(t)-1}\Big(\frac{1}{N}\sum_{i=E^N(t_m)+1}^{E^N(t_{m+1})} \lan \Ir_{C_{a+t-a^N_i}}, \del_{u_i}\ran
- \frac{1}{N}\sum_{t_m<a^N_i\leq t_{m+1}}\oo{\Gs}_c(a+t-a^N_i)
\\
&\quad -\la \Gs_{d, \eps}(a+t-t_{m+1})(t_{m+1}-t_m)\Big)+\frac{\eps}{16}\\
&\leq  \sum_{m=0}^{k(t)-1}\Big( \frac{1}{N}\sum_{i=E^N(t_m)+1}^{E^N(t_{m+1})}\lan \Ir_{C_{a+t-t_{m+1}}}, \del_{u_i}\ran
- \frac{1}{N}\sum_{t_m<a^N_i\leq t_{m+1}}\oo{\Gs}_c(a+t-t_m)
\\
&\quad -\la \Gs_{d, \eps}(a+t-t_m)(t_{m+1}-t_m)\Big)+\frac{\eps}{16}\\
&\leq  \sum_{m=0}^{k(t)-1}\Big( (\bar E^N(t_{m+1})-\bar E^N(t_m))\lan \Ir_{C_{a+t-t_{m+1}}}, \nu^N\ran+\eps_G(N)\\
&\ -(\bar E^N(t_{m+1})-\bar E^N(t_m))\oo{\Gs}_c(a+t-t_m)-\la \Gs_{d, \eps}(a+t-t_{m+1})(t_{m+1}-t_m)\Big)+\frac{\eps}{16},
\end{align*}
where in the first inequality we use \eqref{1000}, for second equality we use \eqref{260} and the fact that
$\oo{\Gs}_c$ is non-increasing,
 in the last inequality we use \eqref{211} and definition of $\Om_G(N)$. By Assumption \ref{cond1}, for any $\eps_1>0$ we have $\rho(\nu^N, \nu)\leq\eps_1$ for all $N$ large.
Therefore for any closed set $B\in\calB([0, \iy))$ and $N$ large
\begin{equation}\label{999}
\nu^N(B)\leq \nu(B^{\eps_1})+\eps_1\ \quad\text{and}\quad \nu(B)\leq \nu^N(B^{\eps_1})+\eps_1.
\end{equation}
Hence on $\Om_A(N)\cap\Om_G(N)$, we have
\begin{align*}
\lan f & , \calG^N_t\ran  - \int_0^t\oo{\Gs}_c(a+t-s)d\bar E^N(s)-\la\int_0^t\oo{\Gs}_d(a+t-s) \, ds\\
& \leq \sum_{m=0}^{k(t)-1}\Big( (\bar E^N(t_{m+1})-\bar E^N(t_m))(\oo{G}(a+t-t_{m+1}-\eps_1)+\eps_1)+\eps_G(N) \\
&\quad -(\bar E^N(t_{m+1})-\bar E^N(t_m))\oo{\Gs}_c(a+t-t_m)-\la \Gs_{d, \eps}(a+t-t_{m+1})(t_{m+1}-t_m)\Big)+\frac{\eps}{16} \\
& \leq \sum_{m=0}^{k(t)-1}\Big( (\bar E^N(t_{m+1})-\bar E^N(t_m))(osc_{\eps_1+\del}(\oo{\Gs}_c)+\eps_1)+\eps_G(N) \\
&\quad + (\bar E^N(t_{m+1})-\bar E^N(t_m))\oo{\Gs}_d(a+t-t_{m+1}-\eps_1)-\la \Gs_{d, \eps}(a+t-t_{m+1})(t_{m+1}-t_m)\Big)+\frac{\eps}{16}\displaybreak
\\
& \leq \sum_{m=0}^{k(t)-1}\Big( (\bar E^N(t_{m+1})-\bar E^N(t_m))(osc_{\eps_1+\del}(\oo{\Gs}_c)+\eps_1)+\eps_G(N) \\
&\quad +\la(t_{m+1}-t_m)\big(\Gs_{d, \eps}(a+t-t_{m+1}-\eps_1)-\Gs_{d, \eps}(a+t-t_{m+1})\big)+2\eps_A(N)|G|_\iy\Big)+\frac{\eps}{8},
\end{align*}
where in the second inequality we use the fact that $\oo{G}=\oo{\Gs}_c+\oo{\Gs}_d$ and in the last inequality we use \eqref{1000} and definition of $\Om_A(N)$.
If we choose $\eps_1<\del$ in \eqref{999}, then
\begin{equation}\label{262}
\sum_{m=0}^{k(t)-1}(t_{m+1}-t_m)(\Gs_{d,\eps}(a+t-t_{m+1}-\eps_1)-\Gs_{d, \eps}(a+t-t_{m+1}))\leq 2N_\eps\del |G|_\iy.
\end{equation}
Hence combining the above estimates we see that for large $N$, on $\Om_A(N)\cap\Om_G(N)$, and $f=\Ir_{C_a}$,
\[
\lann f, \calG^N_t\rann-\int_0^t\oo{\Gs}_c(a+t-s)d\bar E^N(s)-\la\int_0^t\oo{\Gs}_d(a+t-s)ds
\leq  \eps/4,
\]
for all $t\leq T$, if $\del$ is chosen sufficiently small. We can obtain
a lower bound $-\eps/4$ using similar calculations.
The same estimate holds if $\Ir_{C_a}$ is replaced by $\Ir_{\bar C_a}$.
Since $$\lim_{N\to\iy}\p((\Om_A(N))^c\cup(\Om_G(N))^c)=0,$$
 \eqref{202} holds.
Since $\oo{\Gs}_c$ is a bounded continuous function, we have
$$ \sup_{[0, T]}\Big|\int_0^t\oo{\Gs}_c(a+t-s)d\bar E^N(s)-\la\int_0^t\oo{\Gs}_c(a+t-s)ds\Big|\to 0,
\quad\text{in probability.}$$
Now defining
$$H_\la(x, t)= \la\int_0^t\oo{G}(x+t-s)ds, \quad t, x\geq 0,$$
and combining with \eqref{202} we have for any $\eps>0$,
\begin{equation}\label{261}
\lim_{N\to\iy}\p\Big(\sup_{f\in\{\Ir_{C_a}, \Ir_{\bar C_a}\}}\sup_{t\in[0, T]}\Big|\lan f, \calG^N_t\ran-H_\la(a, t)\Big|\geq\eps/2\Big)=0.
\end{equation}
Since $\Gs_{d, \eps}$ has finitely many jumps, there exists $\kap>0$ so that the l.h.s.\ of \eqref{262} is
equal to 0 for all $a\geq\kap$.
Hence the same argument as above gives
\begin{align}\label{1102}
&\lim_{N\to\iy}\p\Big(\sup_{a\geq \kap}\sup_{f\in\{\Ir_{C_a}, \Ir_{\bar C_a}\}}
\sup_{t\in[0, T]}\Big|\lan f, \calG^N_t\ran\notag
\\
&- \int_0^t\oo{\Gs}_c(a+t-s)d\bar E^N(s) -\la\int_0^t\oo{\Gs}_d(a+t-s)ds\Big|\geq\eps/4\Big)=0.
\end{align}
Since $\lim_{x\to\iy}\oo{G}(x)=0$, it is possible to choose $\kap_1, \kap_1\geq \kap\vee\tilde g_{N_\eps},$ such that
$$\oo{G}(x)\leq \frac{\eps}{8(\la T+1)}, \quad \forall x\geq\kap_1. $$
Now for $t\in [0, T]$ and $a\geq \kap_1$, using $\oo{G}=\oo{\Gs}_c+\oo{\Gs}_d$, we get
\begin{align*}
&\Big|\int_0^t\oo{\Gs}_c(a+t-s)d\bar E^N(s)-\la\int_0^t\oo{\Gs}_c(a+t-s)ds\Big|
\\
& =\Big|\int_0^t\big(\oo{G}(a+t-s)-\oo{\Gs}_d(a+t-s)\big) d\bar E^N(s)-\la\int_0^t\big(\oo{G}(a+t-s)-\oo{\Gs}_d(a+t-s)\big)ds\Big|
\\
&\leq \frac{\eps(\bar E^N(T)+\la T)}{8(\la T+1)}+\Big|\la \int_0^t\Big(\oo{\Gs}_d(a+t-s)-\Gs_{d, \eps}(a+t-s))ds
+\Gs_{d, \eps}(a)(\la t-\bar E^N(t))
\\
& \ \, -\int_0^t\big(\oo{\Gs}_d(a+t-s)-\Gs_{d, \eps}(a+t-s)\big)d\bar E^N(s)\Big|
\\
& \leq \frac{\eps(\bar E^N(T)+\la T)}{8(\la T+1)}+\sup_{t\leq T}|\la t-\bar E^N(t)|+\frac{\eps(\bar E^N(T)+\la T)}{16(\la T+1)},
\end{align*}
where in the second inequity we use the fact that $\Gs_{d, \eps}$ is constant on $[\kap_1, \iy)$ and in the last inequality we use \eqref{1000}.
Therefore combining with \eqref{1102}, we have for given $\eps>0$,
\[
\lim_{N\to\iy}\p\Big(\sup_{a\geq \kap_1}\sup_{f\in\{\Ir_{C_a}, \Ir_{\bar C_a}\}}
\sup_{t\in[0, T]}\Big|\lan f, \calG^N_t\ran - H_\la(a,t)\Big|\geq\eps\Big)=0,
\]
for some $\kap_1\geq \kap$. Note that this $\kap_1$ might depend on $\eps$. Hence to complete the proof it is enough to show that
\begin{equation}\label{204}
\lim_{N\to\iy}\p\Big(\sup_{a\in[0, \kap_1]}\sup_{f\in\{\Ir_{C_a}, \Ir_{\bar C_a}\}}
\sup_{t\in[0, T]}\Big|\lan f, \calG^N_t\ran - H_\la(a,t)\Big|\geq\eps\Big)=0.
\end{equation}
Consider $\del>0$ and a partition $\{b_i\}_{i=0}^k$ of $[0, \kap_1]$ so that $\max_{1\leq i\leq k}|b_i-b_{i-1}|\leq\del$. Let $\omega_\del$ be the $\del$-oscillation of $H_\la$ on $[0, \kap_1]\times[0, T]$.
We choose $\del$ small enough so that $\omega_\del<\eps/4$.
Denote
$$\Om_{b_i}(N, \eps)=\Big\{\sup_{f\in\{\Ir_{C_{b_i}}, \Ir_{\bar C_{b_i}}\}}\sup_{t\in[0, T]}\Big|
\lan f, \calG^N_t\ran-H_\la(b_i,t)\Big|\geq\eps\Big\}.
$$
Then from \eqref{261} we have for any $\eps>0$,
\begin{equation}\label{205}
\lim_{N\to\iy}\p(\Om_{b_i}(N, \eps))=0
\end{equation}
for all $i=1,\ldots, k$.
For any $a\in[b_i, b_{i+1}], i=0,\ldots, k-1$, $f\in\{\Ir_{C_a}, \Ir_{\bar C_a}\}$,
\begin{align*}
\sup_{t\in[0, T]}|\lan f, \calG^N_t\ran - H_\la(a, t)|
&\leq\sup_{t\in[0, T]}\Big|\lann \Ir_{\bar C_{b_i}}, \frac{1}{N}\sum_{i=0}^{E^N(t)}
\del_{L_i^N(t)}\rann - H_\la(b_i, t)\Big|\\
& \quad +\sup_{t\in[0, T]}|\lan \Ir_{ C_{b_{i+1}}}, \calG^N_t\ran - H_\la(b_{i+1}, t)|+\omega_\del.
\end{align*}
The above inequality implies that the event in \eqref{204} is contained in
\[
\cup_{i=1}^k(\Om_{b_i}(N,\eps/8)\cup\Om_{b_{i+1}}(N,\eps/8)).
\]
Hence the claim \eqref{204} follows using \eqref{205}.\hfill $\Box$

Define the potential queue-length process
\begin{equation}\label{302}
\calH^N_t(B)=\calQ^N_0(B+t)+\lann \Ir_B,\sum_{i=1}^{E^N(t)}\del_{L_i^N(t)}\rann,\quad B\in\calB([0,\iy)), \quad t\geq 0.
\end{equation}
Then $\calH^N_t(0, \iy)$ represents the total number
of customers that have arrived before time $t\geq 0$ and have positive lead times at
time $t$. The following lemma establishes the law of large number limit for $\calH^N$.
Recall the definition of $H$ in \eqref{301} and that we have set $\oo{G_0}(a)=\calQ_0(a,\iy)$.

\begin{lemma}\label{lem-3}
For any given $\eps>0$, we have
$$\lim_{N\to\iy}\p\Big(\sup_{a\in\R_+}\sup_{f\in\{\Ir_{C_a}, \Ir_{\bar C_a}\}}\sup_{t\in[0, T]}
\Big|\frac{1}{N}\lan f, \calH^N_t\ran-H(a, t)\Big|\geq \eps\Big)=0.$$
\end{lemma}

\noi{\bf Proof:} In view of Lemma \ref{lem-2}, it is enough to show that for any $\eps>0$,
\begin{equation}\label{T0}
\lim_{N\to\iy}\p\Big(\sup_{a\in\R_+}\sup_{f\in\{\Ir_{C_a}, \Ir_{\bar C_a}\}}
\Big|\frac{1}{N}\lan f, \calQ^N_0\ran-\oo{G_0}(a)\Big|\geq \eps\Big)=0.
\end{equation}
By Assumption \ref{initial-meas}, we see that $\rho(\frac{1}{N}\calQ^N, \calQ_0)\to 0$ in probability. Let $\eps_2\in (0, \eps/4)$ be such that
$osc_{\eps_2}(\oo{G_0})<\eps/4$. Then $\rho(\frac{1}{N}\calQ^N, \calQ_0)<\eps_2$ implies that for any $a\in \R_+$ and $f\in\{\Ir_{C_a}, \Ir_{\bar C_a}\}$ we have
\begin{align*}
\frac{1}{N}\lan f, \calQ^N_0\ran-\oo{G_0}(a)& \leq \oo{G_0}(a-\eps_2)-\oo{G_0}(a)+\eps_2\leq osc_{\eps_2}(\oo{G_0})+\eps_2,
\\
\frac{1}{N}\lan f, \calQ^N_0\ran-\oo{G_0}(a)&\geq \oo{G_0}(a+\eps_2)-\oo{G_0}(a)-\eps_2\geq osc_{\eps_2}(\oo{G_0})-\eps_2.
\end{align*}
Hence our required event is a subset of $\{\rho(\frac{1}{N}\calQ^N, \calQ_0)\geq\eps_2\}$.
The proof follows from Assumption \ref{initial-meas}(iii). \hfill $\Box$

The following lemma gives an upper estimate on the frontier process.
\begin{lemma}\label{lem-4}
Let $T, \eps>0$ be given. Define $\Om_{F, \eps}=\{\sup_{t\in[0,T]}F^N(t)\leq y^*+\eps\}$,
where $y^*$ is given in \eqref{y-cond}.
If $\la>\mu$ then
\[
\lim_{N\to\iy}\p(\Om_{F, \eps}^c)=0.
\]
\end{lemma}

\noi{\bf Proof:}\ Take $\del\in(0, \eps/2)$. Then
\begin{align}\label{207}
&\p\Bigl(\sup_{t\in[0,T]}F^N(t)\geq y^*+\eps\Bigr)\notag
\\
&= \p\Bigl(\sup_{t\in[0,\del]}F^N(t)\geq y^*+\eps\Bigr)+\p\Bigl(\sup_{t\in[\del,T]}F^N(t)\geq y^*+\eps, \sup_{t\in[0,\del]}F^N(t)< y^*+\eps\Bigr).
\end{align}
Applying Lemma \ref{lem-3} for $T=0$, specifically, \eqref{T0},
we have for any positive $\eta$,
$$\limsup_{N\to\iy}\p\Bigl(|\frac{1}{N}\calQ^N_0((\del, y^*+\eps/2])-\calQ_0((\del, y^*+\eps/2])|\geq\eta\Bigr)=0. $$
Because of Assumption \ref{initial-meas}, we have $\calQ_0([0, y^*+\eps/2])>0$ and thus we can choose $\del>0$ to satisfy
$2\mu\del<\calQ_0((\del, y^*+\eps/2])$. Thus
\begin{equation}\label{208}
\lim_{N\to\iy}\p\Bigl(\frac{1}{N}\calQ^N_0(\del, y^*+\eps/2])\geq2\mu\del-\eta\Bigr)=1.
\end{equation}
Again Assumption \ref{initial-meas}(iii) implies that $\lim_{N\to\iy}\p(F^N(0)\geq y^*+\eps/2)=0$.
If $F^N$ jumps upward then at the time of jump a customer leaves the queue (due to service/lack of patience) and there is
another customer in the queue waiting to be served.
On $\Om_S(N)\cap\{F^N(0)< y^*+\eps/2\}$,
if $\sup_{t\in[0,\del]}F^N(t)\geq y^*+\eps $ then there exists $\tau^N, 0<\tau^N\leq\del,$ such that
$C^N(\tau^N)>y^*+\eps/2$. Therefore all the customers at time $0$ with their lead times in $(\del, y^*+\eps/2+\tau^N]$ should leave
the system by time $\tau^N$. Again the customers at time $0$ with their initial lead times in $(\del, y^*+\eps/2+\tau^N]$ can leave the
queue by time $\del$ if and only if they receive service. As a result we have
\[
D^N(0,\del]+1\geq \calQ^N_0((\del, y^*+\eps/2+\tau^N]),
\]
hence
\[
\mu\del+\eps_S(N)+\frac{1}{N}\geq
\frac{1}{N}\calQ^N_0((\del, y^*+\eps/2+\tau^N])
\ge\frac{1}{N}\calQ^N_0((\del, y^*+\eps/2]).
\]
Now choosing $\eta=\mu\del/2$, we see, using \eqref{212} and \eqref{208}, that
\begin{equation}\label{209}
\lim_{N\to\iy}\p\Bigl(\sup_{t\in[0,\del]}F^N(t)\geq y^*+\eps\Bigr)=0.
\end{equation}
Fix the above choice of $\del$ and consider the second term on the r.h.s.\ of \eqref{207}.
Now take $\eta<\del$. On $\Om_S(N)\cap\Om_G(N)\cap\Om_A(N)$, for $L=M=\la T+1$, if $\sup_{t\in[\del, T]}F^N(t)\geq y^*+\eps$ and $\sup_{t\in[0, \del]}F^N(t)< y^*+\eps$ then
there is a time
$\tau^N, \del<\tau^N\leq T,$ such that  $C^N(\tau^N)>y^*+\eps/2$. Therefore
 all the customers that arrived in the time interval $(\tau^N-\eta, \tau^N]$ with their lead times at $\tau^N$
in $(0, y^*+\eps/2]$ should leave the queue by time $\tau^N$. As a result we have

\begin{equation*}
D^N(\tau^N-\eta, \tau^N]+1\geq \sum_{i=E^n(\tau^N-\eta)+1}^{E^N(\tau^N)}\del_{L_i^N(\tau^N)}((0, y^*+\eps/2])
\geq \sum_{i=E^n(\tau^N-\eta)+1}^{E^N(\tau^N)}\del_{u_i^N}((\eta, y^*+\eps/2]),
\end{equation*}
and therefore
\[
\mu\eta+2\eps_S(N)+\frac{1}{N}\geq \eta\la\nu^N(\eta, y^*+\eps/2]-2\eps_A(N)-2\eps_G(N).
\]
Since $\rho(\nu^N, \nu)\to 0$, we have $\nu^N(\eta, y^*+\eps/2]\geq
\nu[3\eta/2, y^*+\eps/2-\eta/2]-\eta/2$ for large $N$.
Now $G$ being strictly increasing
 we can choose $\eta$ small enough so that $\la\nu[3\eta/2, y^*+\eps/2-\eta/2]-\la\eta/2-\mu>0$. Thus for large $N$,
\begin{equation}\label{T1}
\Bigl\{\sup_{t\in[\del, T]}F^N(t)\geq y^*+\eps,\sup_{t\in[0, \del]}F^N(t)\geq y^*+\eps\Bigr\}\subset (\Om_S(N)\cap\Om_G(N)\cap\Om_A(N))^c.
\end{equation}
Hence $\lim_{N\to\iy}\p(\sup_{t\in[\del,T]}F^N(t)< y^*+\eps, \sup_{t\in[0,\del]}F^N(t)< y^*+\eps)=0$. The result follows
now from \eqref{207}, \eqref{209} and \eqref{T1}. \hfill $\Box$

The following lemma shows that in fluid limit the number of customers in the queue at time $t$ with their lead times in $[C^N(t), F^N(t))$ is negligible. The proof of the lemma is in the spirit of \cite[Proposition 3.6]{Doy-Leh-Shre} where a similar result
is obtained for a single-server queueing model under preemptive EDF-e policy.

\begin{lemma}\label{lem-5}
Let $\eps>0$ be given. Then for $T>0$,
\begin{equation}\label{250}
\lim_{N\to\iy}\p\Bigl(\sup_{t\in[0,T]}\frac{1}{N}\calQ^N_t[C^N(t),F^N(t)]>\eps\Bigr)=0.
\end{equation}
\end{lemma}

\noi{\bf Proof:}  First we show that
\begin{equation}\label{r02}
\lim_{N\to\iy}\p\Bigl(\sup_{t\in[0,T]}\frac{1}{N}\calQ^N_t[C^N(t),F^N(t))>\eps\Bigr)=0.
\end{equation}
Define
$$ \sig^N(t)=\sup\{s\in[0, t]\ :\ C^N(s)=F^N(s)\}.$$
Consider $t>0$ so that $\frac{1}{N}\calQ^N_t[C^N(t),F^N(t))>\eps$.
First we note that
\[
\calQ^N_{\sig^N(t)}([C^N(\sig^N(t)), F^N(\sig^N(t)))=1.
\]
This is true because both $C^N, F^N$ have RCLL paths that decay at unit rate,
and therefore there must be an arrival at
time $\sig^N(t)$. Since $C^N<F^N$ on $(\sig^N(t), t]$, $F^N(s)= F^N(\sig^N(t))+s-\sig^N(t), \ s\in [\sig^N(t), t]$.
Therefore the number of customers in the queue at time $t$ with their lead times in $[C^N(t),F^N(t))$ who were also in queue
at time $\sig^N(t)$ must not be bigger than $\calQ^N_{\sig^N(t)}[C^N(\sig^N(t)), F^N(\sig^N(t)))$. Thus
\begin{equation}\label{213}
\calQ^N_t[C^N(t),F^N(t))\leq 1+\sum_{i=E^N(\sig^N(t))+1}^{E^N(t)}\del_{L^i(t)}((0, F^N(t)))-D^N(\sig^N(t), t].
\end{equation}
First we note that $(\sig^N(t), t]$ is included in the busy period starting at $\sig^N(t)$. Therefore
\begin{align}\label{214}
\frac{1}{N}D^N(\sig^N(t), t] &= \bar S^N\Bigl(\int_0^t\Ir_{\bar X^N(s)>0}ds\Bigr)-\mu\int_0^t\Ir_{\bar X^N(s)>0} ds\notag
\\
&\ -\bar S^N\Bigl(\int_0^{\sig^N(t)}\Ir_{\bar X^N(s)>0}ds\Bigr)+\mu\int_0^{\sig^N(t)}\Ir_{\bar X^N(s)>0}ds-\mu(t-\sig^N(t)).
\end{align}
On $\Om_{F,\del}(N)$
(for $\la\leq\mu$, take $y^*=\iy$), for some $\del>0$,
\begin{align}\label{215}
\sum_{i=E^N(\sig^N(t))+1}^{E^N(t)}\del_{L_i^N(t)}((0, F^N(t)))&\leq \sum_{i=E^N(\sig^N(t))+1}^{E^N(t)}\del_{u_i^N}((0, F^N(t)+t-\sig^N(t)))\notag
\\
& =\sum_{i=E^N(\sig^N(t))+1}^{E^N(t)}\del_{u_i^N}((0, F^N(\sig^N(t))))\notag
\\
& \leq \sum_{i=E^N(\sig^N(t))+1}^{E^N(t)}\del_{u_i^N}((0, y^*+\del)).
\end{align}
Hence using \eqref{213}, \eqref{214} and \eqref{215}, on $\Om_S(N)\cap\Om_G(N)\cap\Om_A(N)\cap\Om_{F, \del}$, for large $N$
and $\eps_1$ as in \eqref{999}, we have
\begin{align*}
\frac{1}{N}\calQ^N_t[C^N(t),F^N(t))&\leq
\frac{1}{N}+ \la (t-\sig^N(t))\nu^N[0, y^*+\del]-\mu(t-\sig^N(t))\\
&\quad +2(\eps_A(N)+\eps_G(N)+\eps_S(N))
\\
&\leq\frac{1}{N}+ \la (t-\sig^N(t))\nu[0, y^*+\del+\eps_1]+\la T\eps_1-\mu(t-\sig^N(t))
\\
& \quad +2(\eps_A(N)+\eps_G(N)+\eps_S(N))
\\
&\leq \frac{1}{N}+\lambda T\nu(y^*, y^*+\del+\eps_1]+\la T\eps_1+2(\eps_A(N)+\eps_G(N)+\eps_S(N)),
\end{align*}
where in the second inequality we used Assumption \ref{cond1}(iii) and in the
third inequality the fact that $\la\nu[0, y^*]=\mu$. Note that $\eps_1>0$ can be chosen arbitrarily small.
Hence taking $\del, \eps_1>0$ small and using Assumption \ref{G-cond}, we see that for large $N$,
$$ \Bigl\{\sup_{t\in[0,T]}\frac{1}{N}\calQ^N_t[C^N(t),F^N(t))>\eps\Bigr\}\subset (\Om_S(N)\cap\Om_G(N)\cap\Om_A(N)\cap\Om_{F, \del})^c.$$
Hence \eqref{r02} follows using \eqref{201},\eqref{212},
\eqref{90} and Lemma \ref{lem-4}.

Now to complete the proof
it is enough to see that the following holds:
\[
\sup_{[0, T]}\frac{1}{N}\calQ^N_t(\{F^N(t)\})\leq \sup_{[0, T]}\frac{1}{N}\calH^N_t(\{F^N(t)\})\to 0\ \text{in probability},
\]
as $N\to\iy$,where we use Lemma \ref{lem-3} and Assumption \ref{initial-meas}(i)
for the last claim. \hfill $\Box$

The following lemma plays a key role in our analysis. In Lemma \ref{lem-5} we have seen that the lead time
tends to coincide with the frontier asymptotically. Therefore it is likely that there is hardly any reneging
when the frontier is away from zero. The following lemma establishes this fact.

\begin{lemma}\label{lem-5.5}
Let $T>0$ be given. Then for any $\eps, \del>0$ ,
we have
\[
\lim_{N\to\iy}\p\Bigl(\frac{1}{N}\int_0^T\Ir_{\{F^N(s)>\del\}}dR^N(s)\geq\eps\Bigr)=0.
\]
\end{lemma}

\noi{\bf Proof:}\ We define excursion times as follows:
\begin{align*}
\sig^N_1 &= \inf\{t\geq 0\ :\ F^N(s)\geq\del\}\wedge T, \\
 \tilde\sig^N_i &= \inf\{s\geq\sig^N_i\ :\ F^N(s)\leq \del/2 \}\wedge T,\\
\sig^N_i &= \inf\{s\geq\tilde\sig^N_{i-1} \ :\ F^N(s)\geq \del \}\wedge T, \ i\geq 2.
\end{align*}

Recall $H(\cdot, \cdot), \calH^N$ from \eqref{301}, \eqref{302}. From
 Lemma \ref{lem-3}, we obtain a
sequence $\{\eps_H(N)\}$ with $\eps_H(N)\to 0\ \mbox{as}\ N\to\iy,$ such that
\begin{equation}\label{777}
\p\Bigl(\sup_{a\in\R_+}\sup_{f\in\{\Ir_{C_a}, \Ir_{\bar C_a}\}}\sup_{t\in[0, T]}\Bigl|\frac{1}{N}\lan f,
\calH^N_t\ran-H(a, t)\Bigr|\leq\eps_H(N)\Bigr)\geq 1-\eps_H(N).
\end{equation}
Denote
\begin{equation}\label{91}
\Om_H(N)=\{\sup_{a\in\R_+}\sup_{f\in\{\Ir_{C_a}, \Ir_{\bar C_a}\}}\sup_{t\in[0, T]}|\frac{1}{N}\lan f, \calH^N_t\ran-H(a, t)|\leq\eps_H(N)\}.
\end{equation}
We divide the proof to the following two cases.

\noindent{\bf Case 1:}\ $\sig^N_1=T$.
On $\Om_H(N)$, since $\sig^N_1=T$,
\begin{align*}
\frac{1}{N}\int_0^T \Ir_{\{F^N(s)>\del\}}dR^N(s)&\leq
\frac{1}{N}(R^N(T)-R^N(T-))
\\
&\leq
\frac{1}{N}\calQ^N_0(\{T\})+\frac{1}{N}\sum_{i=1}^{E^N(T)}
\Ir_{L_i^N(T)}\{0\}
\\
&=\frac{1}{N}\calH^N_T(\{0\})
\\
&= \frac{1}{N}\lan \Ir_{\bar C_0}, \calH^N_T\ran-\frac{1}{N}\lan\Ir_{ C_0}, \calH^N_T\ran \leq 2\eps_H(N),
\end{align*}
where the second inequality follows from the continuity at $0$ of $a\mapsto H(a,T)$
(see the discussion after Remark \ref{rem1}).
Hence for all $N$ large, on $\Om_H(N)$, we have
\begin{equation}\label{888}
\frac{1}{N}\int_{0}^T\Ir_{\{F^N(s)>\del\}}dR^N(s)\leq 2\eps_H(N)< \frac{\eps}{2}.
\end{equation}

\noindent{\bf Case 2:}\ $\sig^N_1<T$.
Note that $F^N(\cdot)$ decreases continuously at a unit rate in between successive jumps,
and always jumps upwards. We therefore have $(\tilde\sig^N_i-\sig^N_i)\geq \del/2$
and for all $i$, $F^N(\tilde\sig^N_i)=\del/2$, $F^N(\sig^N_i)\geq\del$
and finally, $F^N(s)< \del$ for $s\in[\tilde\sig^N_i, \sig^N_{i+1})$. This implies that
$i$ runs over a finite set of integers. In fact, $\max\{i: \sig^N_i\leq T\}\leq \frac{2T}{\del}$.
Also
$\Ir_{\{F^N(\sig^N_i)>\del\}}(R^N(\sig^N_i)-R^N(\sig^N_i-))$ is positive only if at time $\sig^N_i$ one of the customers abandons the system
 and $F^N(\sig^N_i)>\del$ holds. We shall show that this can happen only on a set of negligible measure.
Now
\begin{align}\label{217}
\frac{1}{N}\int_{0}^T\Ir_{\{F^N(s)>\del\}}dR^N(s) &\leq \frac{1}{N}\sum_{i}\int_{[\sig^N_i,\tilde\sig^N_i]}dR^N(s)\notag
\\
&=\sum_{i}\frac{1}{N}(R^N(\sig^N_i)-R^N(\sig^N_i-))+\sum_{i}\frac{1}{N}(R^N(\tilde\sig^N_i)-R^N(\sig^N_i)).
\end{align}

Let $J(\del)$ be the maximum number of intervals of the form $[\sig^N_i, \tilde\sig^N_i]$ that are subsets of $[0 , T]$. We have shown above that $J(\del)\leq \frac{2T}{\del}$. Take $\kap<\del/2$ positive. We subdivide each $[\sig^N_i, \tilde\sig^N_i]$
into intervals of length $\kap$ with at most one subinterval of length less than $\kap$. We denote these intervals by $[\sig^N_{i, j},\sig^N_{i, j+1}]$ where $j$ varies over suitable number of indices. Let $J_1(\kap)$ denote the sum over $i$ of the number of intervals of this
form. We note that $J_1(\kap)\leq \frac{2T}{\del}(\lfl\frac{\del}{2\kap}\rfl+1)$ . Hence from \eqref{217}, we obtain
\begin{equation}\label{218}
\frac{1}{N}\int_{0}^T\Ir_{\{F^N(s)>\del\}}dR^N(s)\leq
\sum_{i}\frac{1}{N}(R^N(\sig^N_i)-R^N(\sig^N_i-))+\sum_{i, j}\frac{1}{N}(R^N(\sig^N_{i, j+1})-R^N(\sig^N_{i, j})).
\end{equation}
Now for any fixed $i, j$ we obtain
\begin{equation*}
\frac{1}{N}\Bigl(R^N(\sig^N_{i, j+1})-R^N(\sig^N_{i, j})\Bigr)
\leq \frac{1}{N}\calQ^N_{\sig^N_{i, j}}[C^N(\sig^N_{i, j}), F^N(\sig^N_{i, j}))
+ \frac{1}{N}\sum_{i=E^N(\sig^N_{i, j})}^{E^N(\sig^N_{i, j+1})}\Ir_{L_i^N(\sig^N_{i, j+1})}(-\kap, 0].
\end{equation*}
To see that the above holds it is enough to count the number of customers who can loose their patience in the time interval $(\sig^N_{i, j},\sig^N_{i, j+1}]$. We have $(\sig^N_{i, j+1}-\sig^N_{i, j})\leq\kap< F^N(\sig^N_{i, j})$ as $F^N\geq\del/2$ on $[\sig^N_i, \tilde\sig^N_i]$. Therefore
customers who were in the queue at time $\sig^N_{i, j}$ with their lead times greater than $ F^N(\sig^N_{i, j})$ can not loose their
patience in the interval $(\sig^N_{i,j},\sig^N_{i, j+1}]$. The rightmost term above counts the number of possible new customers who can renege
the system. Hence for large $N$, on $\Om_A(N)\cap\Om_G(N)$,
\begin{align*}
\frac{1}{N}(R^N(\sig^N_{i, j+1})-R^N(\sig^N_{i, j}))&\leq
\sup_{[0, T]}\frac{1}{N}\calQ^N_t[C^N(t), F^N(t))
\\
&\quad +\la(\sig^N_{i, j+1}-\sig^N_{i, j})\nu^N[0, \kap] +2(\eps_{A}(N)+
\eps_G(N))
\\
&\leq\sup_{[0, T]}\frac{1}{N}\calQ^N_t[C^N(t), F^N(t))
\\
&\quad +\la(\sig^N_{i, j+1}-\sig^N_{i, j})(\nu[0, \kap+\eps_1]+\eps_1) +2(\eps_{A}(N)+
\eps_G(N)),
\end{align*}
where in the last line we use \eqref{999} for some $\eps_1>0$.
Hence from \eqref{218} we have, on $\Om_A(N)\cap\Om_G(N)\cap\Om_H(N)$,
\begin{align}\label{219}
\frac{1}{N}\int_0^T\Ir_{\{F^N(s)>\del\}}dR^N(s)&\leq 2J(\del)\eps_H(N)+ J_1(\kap)\Big[\sup_{[0, T]}\frac{1}{N}\calQ^N_t[C^N(t), F^N(t))+2(\eps_{A}(N)\notag
\\
&\ +\eps_G(N))\Big]+\la T\nu[0, \kap+\eps_1]+\la T\eps_1,
\end{align}
where the first sum in \eqref{218} is estimated as in case 1.
We next use
Lemma \ref{lem-5} to get a sequence $\{\eps_{CF}(N)\}, \eps_{CF}(N)\to 0 \ \mbox{as}\ N\to\iy,$ satisfying
$$ \p\Bigl(\sup_{[0, T]}\frac{1}{N}\calQ^N_t[C^N(t), F^N(t)]\leq\eps_{CF}(N)\Bigr)\geq 1-\eps_{CF}(N).$$
Denote
\begin{equation}\label{92}
\Om_{CF}(N)=\{\sup_{[0, T]}\frac{1}{N}\calQ^N_t[C^N(t), F^N(t))\leq\eps_{CF}(N) \}.
\end{equation}
To conclude we now choose first $\kappa, \eps_1$ small enough so as to match the last two terms from \eqref{219} and then $N$ large so that on
$\Om_A(N)\cap\Om_G(N)\cap\Om_H(N)\cap\Om_{CF}(N)$, from \eqref{219},
\begin{equation}\label{220}
\frac{1}{N}\int_{0}^T\Ir_{\{F^N(s)>\del\}}dR^N(s)\leq \frac{\eps}{4}.
\end{equation}
Hence the proof follows from \eqref{888} and \eqref{220} by observing that $\p(\Om_H(N)\cap\Om_A(N)\cap\Om_G(N)\cap\Om_{CF}(N))\to 1$ as $N\to\iy$. \hfill $\Box$

The next result shows that under Assumption \ref{initial-meas}, the system will remain busy with high probability if $\la\geq\mu$.

\begin{lemma}\label{lem-6}
Let Assumption \ref{initial-meas} hold. Fix $T>0$. Then
for $\la\geq\mu$ we have
\[
\liminf_{N\to\iy}\p(\inf_{t\in[0, T]}\frac{1}{N}Q^N(t)>0)=1.
\]
\end{lemma}

\noi{\bf Proof:}\ Since $\frac{1}{N}Q^N(0)\Ra Q(0)$ and $Q(0)>0$, we have $\liminf_{N\to\iy}\p(\frac{1}{N}Q^N(0)\geq\Del)=1$
where $2\Del=Q(0)$. Denote
\begin{equation}\label{94}
\Om_\Delta(N)=\{\frac{1}{N}Q^N(0)\geq\Del\}.
\end{equation}
Choose $0<\del<\frac{\Del}{2}$. Define
$$\bar\tau^N=\inf\Bigl\{t\geq 0\ :\ \frac{1}{N}Q^N(t)=0\Bigr\}\ \mbox{and}\ \bar\sig^N=\sup\Bigl\{t\leq\bar\tau^N\wedge T\ :\ \frac{1}{N}Q^N(t)\geq\del\Bigr\}. $$
By definition $\frac{1}{N}Q^N(t)\leq \del$ for $t\in(\bar\sig^N, \bar\tau^N\wedge T]$.
It is easy to see that $0<\bar\sig^N\leq\bar\tau^N\wedge T$, on $\Om_\Delta(N)$.
Recall $\Om_H(N)$ from \eqref{777}. Then
$$|\frac{1}{N}Q(\bar\sig^N)-\frac{1}{N}Q(\bar\sig^N-)|\leq \frac{1}{N}\Big(1+\calQ^N_{\bar\sig^N-}(\{0\})\Big), $$
where the r.h.s.\ denotes maximum number of customers who can leave at time $\bar\sig^N$. Therefore on $\Om_H(N)$, for
large $N$, we have
$$|\frac{1}{N}Q(\bar\sig^N)-\frac{1}{N}Q(\bar\sig^N-)|\leq \frac{1}{N}\Big(1+\calH^N_{\bar\sig^N-}(\{0\})\Big)\leq \frac{1}{N}+2\eps_H(N)
< \frac{\del}{2}. $$
Hence for $N$ large, on $\Om_\Delta(N)\cap\Om_H(N)\cap\{\bar\tau^N\leq T\}$, we have $0<\bar\sig^N<\bar\tau^N$ and $\frac{1}{N}Q^N(\bar\sig^N)\geq\frac{\del}{2}$.
 Now we prove that if $\min_{[0, t]}Q^N(s)>0$, then for $a\geq 0$,
\begin{equation}\label{786}
\calQ^N_t((F^N(t)\vee a, \iy))= \calH^N_t((F^N(t)\vee a, \iy)).
\end{equation}
It is easy to see that the r.h.s.\ is bigger than the l.h.s. For the other inequality, consider any customer (say $i$-th)
 that has arrived before time $t$, with lead time
at time $t$, that is greater
that $F^N(t)\vee a$. We show that the customer is still in queue. Since $F^N(t)\geq F^N(s)-t+s$ for all $0\leq s\leq t$, by
definition of $F^N(\cdot)$ the $i$-th customer has never become a head-in-the-line customer in $[0, t]$. Hence the possible
ways in which this $i$-th customer may not be present in the l.h.s.\ are as follows: (a) at some time $s\leq t$, the server finishes a job and
the $i$-th customer arrives with deadline less than $C^N(s-)$ and directly goes into service, (b) at some time $s\leq t$, the server finishes a job and a customer reneges the queue. If the $i$-th customer is next to the head-in-the-line, then he goes into service without
even being the head-in-line. In both cases, it is easy to see that $C^N(s)$ is bigger than the lead time of customer-$i$ at the time $s$.
Thus $F^N(s)$ is bigger than the lead time of the $i-$th customer at time $s$. But this is contradicting the fact that the lead time
of customer-$i$ at time $t$ is bigger than $F^N(t)$. This proves claim \eqref{786}.

Hence for any $t<\bar\tau^N$, $Q^N(t)>0\Ra F^N(t)>0,$ and (by \eqref{786})
$$\calQ^N_t((F^N(t), \infty))=\calH^N_t((F^N(t), \infty))=\calQ^N_0((t+F^N(t), \infty))+\sum_{i=1}^{E^N(t)}\Ir_{L^i(t)}(F^N(t), \iy).$$
Note that $\min_{[0, T]}H(0, t)>2\eps$ for some $\eps>0$.
Therefore on $\Om_H(N)$, for large $N$, \[\frac{1}{N}\calH^N(t)((0,\iy))\geq\eps.\]
Choose $\del>0$ small enough so that $\del<\eps/2$. Hence on $\Om_H(N)\cap\Om_\Delta(N)$, for $N$ large and $t\in[\bar\sig^N, \bar\tau^N\wedge T)$,
\begin{align*}
 \frac{1}{N}Q^N(t)\leq \del <  \frac{1}{N}\calH^N(t)((0, \iy))-\del
&=\frac{1}{N}\calH^N_t((0, F^N(t)])+\frac{1}{N}\calH^N_t((F^N(t), \iy))-\del
\\
&=\frac{1}{N}\calH^N_t((0, F^N(t)])+\frac{1}{N}\calQ^N_t((F^N(t), \iy))
-\del.
\end{align*}
Hence
\[
\frac{1}{N}\calH^N_t((0, F^N(t)])>\del,
\]
and thus
\[
2\eps_H(N)+\calQ_0(t, F^N(t)+t]+ \la\int_0^t\nu(t-s, F^N(t)+t-s]ds\geq \del.
\]
Using the continuity of $H$ we see that $\inf\{x\in\R_+\ : \ \inf_{[0, T]}(\calQ_0(t, x+t]+ \la\int_0^t\nu(t-s, x+t-s]ds)\geq \del/2\}>0$. Therefore we can choose $\del_1\in(0, \del/2)$, depending on the lower bound of the above infimum, such that
$$\calQ_0(t, x+t]+ \la\int_0^t\nu(t-s, x+t-s]ds\geq \del/2\quad \text{implies}\quad x> \del_1,$$
and thus for all $N$ large, on $\Om_H(N)\cap\Om_\Delta(N)$, for $t\in [\bar\sig^N, \bar\tau^N\wedge T)$,
$$\frac{1}{N}Q^N(t)\leq \del\Ra F^N(t)>\del_1.$$
Again by definition of $F^N$, we see that $F^N(t)>\del_1$ on $[\bar\sig^N, \bar\tau^N\wedge T)$ implies that $F^N(\bar\tau^N\wedge T)\geq \del_1$.
Therefore on $\Om_H(N)\cap\Om_\Delta(N)$, for $N$ large, we have
\begin{equation}\label{216}
R^N(\bar\tau^N\wedge T)-R^N(\bar\sig^N)=\int_{\bar\sig^N}^{\bar\tau^N\wedge T}\Ir_{\{\frac{1}{N}Q^N(s)\leq\del\}}dR^N(s)
\leq \int_{\bar\sig^N}^{\bar\tau^N\wedge T}\Ir_{\{F^N(s)>\del_1\}}dR^N(s).
\end{equation}

We denote the complement of the event in Lemma \ref{lem-5.5} by $\Om_{RF}(N)$ for $\eps=\frac{\del}{4}, \del=\del_1$.
Then $\p(\Om_{RF}(N))\to 1$ as $N\to\iy$ and on $\Om_{RF}(N)\cap\Om_H(N)\cap\Om_\Delta(N)$,
\begin{equation}\label{221}
\frac{1}{N}(R^N(\bar\tau^N\wedge T)-R^N(\bar\sig^N))\leq \frac{\del}{4}.
\end{equation}
Now we consider the equations on $[\bar\sig^N, \bar\tau^N\wedge T]$ for large $N$, on the event
$\Om_{RF}(N)\cap\Om_A(N)\cap\Om_S(N)\cap\Om_H(N)\cap\Om_\Delta(N)\cap\{\bar\tau^N\leq T\}$, i.e.,

\begin{align*}
0=\frac{1}{N}Q^N(\bar\tau^N)&=\frac{1}{N}Q^N(\bar\sig^N)+\frac{1}{N}E^N(\bar\sig^N, \bar\tau^N]-\frac{1}{N}D^N(\bar\sig^N, \bar\tau^N]
-\frac{1}{N}(R^N(\bar\tau^N)-R^N(\bar\sig^N))
\\
&\geq \frac{\del}{2}-2\eps_A(N)+\la(\bar\tau^N-\bar\sig^N)-2\eps_S(N)-\mu(\bar\tau^N-\bar\sig^N)-\frac{\del}{4}
\\
&\geq \frac{\del}{4}-2(\eps_A(N)+\eps_S(N)),
\end{align*}
where in the last line we have used the fact $\la\geq\mu$ . But this is a contradiction for all $N$ large enough.
Hence our required event $\{\bar\tau^N\leq T\}$ lies in
$\Big(\Om_T(N)\cap\Om_A(N)\cap\Om_S(N)\cap\Om_H(N)\cap\Om_\Delta(N)\Big)^c$.
Hence the proof. \hfill $\Box$

\begin{corollary}\label{cor-1}
Let $\la\geq \mu$. For any $T, \eps>0$,
$$\lim_{N\to\iy}\p\Bigl(\sup_{a\in\R_+}\sup_{[0, T]}\Bigl|\frac{1}{N}\calQ^N_t((F^N(t)\vee a, \iy))-H(F^N(t)\vee a,t)\Bigr|\geq\eps\Bigr)=0.$$
\end{corollary}

\noi{\bf Proof:} If $\inf_{[0, T]}\frac{1}{N}Q^N(t)>0$, then $\frac{1}{N}\calQ^N_t((F^N(t)\vee a, \iy))=\frac{1}{N}\calH^N_t((F^N(t)\vee a, \iy))$\eqref{786}.
Hence the proof follows using Lemma \ref{lem-3} and Lemma \ref{lem-6}.\hfill $\Box$

Before we go to prove the convergence results for the queue length, let us define, for given $\del>0$,
$$ R^{N,\del}_1(t)=\int_0^t\Ir_{\{\frac{1}{N}Q^N(s)<\frac{1}{N}H^N(0,s)-\del\}}dR^N(s)\ \text{and}\
 R^{N,\del}_2(t)=\int_0^t\Ir_{\{\frac{1}{N}Q^N(s)\geq \frac{1}{N} H^N(0,s)-\del\}}dR^N(s),$$
where $H^N(a, s)=\calH^N_s((a,\iy))$ for $a\geq 0$. Therefore for all $\del>0$, $R^N=R^{N,\del}_1+R^{N,\del}_2$.

\begin{lemma}\label{lem-7}
Let $\la\geq\mu$.
For any given $\eps, \del>0$, we have $\lim_{N\to\iy}\p(\sup_{[0, T]}\frac{1}{N}R^{N,\del}_1(t)\geq \eps)=0$.
\end{lemma}

\noi{\bf Proof:}\ Using the same argument as in Lemma \ref{lem-6}, we see that on $\{\inf_{[0, T]}\frac{1}{N}Q^N(t)>0\}\cap
\Om_H(N)$, for all $N$ large,
\[
\frac{1}{N}Q^N(s)<\frac{1}{N}H^N(0,s)-\del\ \mbox{implies}\ F^N(s)\geq\del_1,
\]
for $s\in[0, T]$ and some $\del_1>0$. Hence the proof follows using Lemma \ref{lem-5.5} and Lemma \ref{lem-6}.\hfill $\Box$

\subsection{Proofs for the case $\lambda\geq\mu$}

Now we give a proof of Theorem \ref{main2-th}.

\noi{\bf Proof of Theorem \ref{main2-th}:}
It is enough to show that for any positive $T, \eps,$
\begin{align}
\lim_{N\to\iy}\p\Bigl(\sup_{[0, T]}\Bigl|\frac{1}{N}Q^N(t)-\phi(t)\Bigr|\geq\eps\Bigr) &=0,\label{222}
\\
\lim_{N\to\iy}\p\Bigl(\sup_{[0, T]}\Bigl|\frac{1}{N}R^N(t)-\eta(t)\Bigr|\geq\eps\Bigr) &=0.\label{1004}
\end{align}
Denote $\Om_Q(N)=\{\inf_{[0, T]}\frac{1}{N}Q^N(t)>0\}$.
We know from Lemma \ref{lem-6} that
\[
\lim_{N\to\iy}\p(\Om_Q(N))=1.
\]
Hence from \eqref{102} and \eqref{103} we have on $\Om_Q(N)$ that for any $\del>0$,
\[
\frac{1}{N}Q^N(t)=\frac{1}{N}Q^N(0)+\bar E^N(t)-\bar S^N(t)-\frac{1}{N}R^{N,\del}_1(t)-\frac{1}{N}R^{N,\del}_2(t),
\]
and hence
\begin{align*}
\frac{1}{N}Q^N(t)-\frac{1}{N}H^N(0, t)+\del &=\frac{1}{N}Q^N(0)-\frac{1}{N}H^N(0, t)+\del+\bar E^N(t)-\bar S^N(t)\\
& -\frac{1}{N}R^{N,\del}_1(t)-\frac{1}{N}R^{N,\del}_2(t),
\end{align*}
implying
\begin{align}
\label{223}
-\Bigl(\frac{1}{N}Q^N(t)-\frac{1}{N}H^N(0, t)+\del\Bigr)^-
& = -\Bigl(\frac{1}{N}Q^N(t)-\frac{1}{N}H^N(0, t)+\del\Bigr)^+ +\frac{1}{N}Q^N(0) \notag\\
& -\frac{1}{N}H^N(0, t)+\del +\bar E^N(t)-\bar S^N(t) \notag\\
& -\frac{1}{N}R^{N,\del}_1(t)-\frac{1}{N}R^{N,\del}_2(t).
\end{align}
Define
\begin{align}\label{1002}
Y^{N,\del}(t) &= -\Bigl(\frac{1}{N}Q^N(t)-\frac{1}{N}H^N(0, t)+\del\Bigr)^-,\notag
\\
X^{N,\del}(t) &= -\frac{1}{N}H^N(0, t)+\del-\Bigl(\frac{1}{N}Q^N(t)-\frac{1}{N}H^N(0, t)+\del\Bigr)^+ \notag
 \\
 & \ \ +\frac{1}{N}Q^N(0)
 +\bar E^N(t)-\bar S^N(t)-\frac{1}{N}R^{N,\del}_1(t) .
\end{align}
Hence from \eqref{223} we get
\begin{equation}\label{224}
Y^{N, \del}=X^{N, \del}-\frac{1}{N}R^{N,\del}_2(t),
\end{equation}
on $[0, T]$. We note that $Y^{N, \del}\leq 0$ and $\int_0^\cdot\Ir_{\{Y^{N, \del}(s)<0\}}dR_2^{N,\del}(s)=0$. Hence $Y^{N, \del}$
 is the Skorohod reflection term
 of $X^{N, \del}$ in $(-\iy, 0]$ (see Definition \ref{def1}).
 Therefore using the Lipschitz property
 of the Skorohod map $\Gamma_h, h\equiv 0,$ with Lipschitz constant $2$ we have that
\begin{align}\label{225}
\lefteqn{\sup_{[0, T]}|Y^{N, \del}(t)-\phi(t)+H(0, t)| \leq 2\sup_{[0, T]}|X^{N,\del}-\psi(t)+H(0,t)|}\notag\\
&\leq 2\Big[\sup_{[0, T]}\Bigl|\frac{1}{N}Q^N(0)+\bar E^N(t)-\bar S^N(t)-\frac{1}{N}H^N(0, t)-\psi(t)\notag\\
& +H(0,t)\Bigr|\Big]+\del+\sup_{[0, T]}\Bigl(\frac{1}{N}R^{N,\del}_1(t)+\Bigl(\frac{1}{N}Q^N(t)-\frac{1}{N}H^N(0, t)+\del\Bigr)^+\Bigr)\notag\\
& \leq 2\Big[\sup_{[0, T]}\Bigl|\frac{1}{N}Q^N(0)+\bar E^N(t)-\bar S^N(t)-\frac{1}{N}H^N(0, t)-\psi(t)\notag\\
& +H(0,t)\Bigr|\Big]+\sup_{[0, T]}\frac{1}{N}R^{N,\del}_1(t)+2\del,
\end{align}
where in the second inequality we use \eqref{1002} and for the third inequality we use the fact that $Q^N(t)\leq H^N(0, t)$. Using \eqref{225},
on $\Om_A(N)\cap\Om_S(N)\cap\Om_H(N)$, and the fact that $\psi(t)=Q(0)+(\la-\mu)t$ we get
\begin{align*}
\sup_{[0, T]}\Bigl|\frac{1}{N}Q^N(t)-\phi(t))\Bigr| &\leq\sup_{[0, T]}|Y^{N, \del}(t)-\phi(t)+H(0, t)|+2\del+\sup_{[0, T]}\Bigl|\frac{1}{N}H^N(0,t)-H(0, t)\Bigr|
\\
&\leq 2(\eps_A(N)+\eps_S(N)+2\eps_H(N))+4\del+\Bigl|\frac{1}{N}Q^N(0)-Q(0)\Bigr|\\
&+\sup_{[0, T]}\frac{1}{N}R^{N,\del}_1(t).
\end{align*}
Hence \eqref{222} follows by choosing suitable small $\del>0$ and applying Assumption \ref{initial-meas}(iii) and Lemma \ref{lem-7}.

To prove \eqref{1004} we observe that by \eqref{102}, \eqref{103} and \eqref{105}
one has
$$\sup_{[0, T]}|\frac{1}{N}R^N-\eta|\leq \sup_{[0, T]}|\frac{1}{N}Q^N-\phi|+ \sup_{[0, T]}|\frac{1}{N}Q^N(0)+\bar E^N-D^N-\psi|+\frac{2}{N}, $$
and the middle term on r.h.s.\ goes to zero using the fact that $D^N=S^N$ on $\{\inf_{s\in[0, T}Q^N(s)>0\}$ and using Lemma \ref{lem-6}. \hfill $\Box$

\skp

\noi{\bf Proof of Theorem \ref{main-th}(a):}\ Since $Q^N(t)= \calQ^N_t([C^N(t), F^N(t)])+\calQ^N_t((F^N(t), \iy))$, using \eqref{250}
and Corollary \ref{cor-1}, we see that $\frac{1}{N}Q^N(\cdot)-H(F^N(\cdot), \cdot)\Ra 0$ as $N\to\iy$. Hence applying Theorem \ref{main2-th} we obtain that $H(F^N(\cdot), \cdot)\Ra\phi(\cdot)$ as $N\to\iy$. At this point we note that $\phi(t)\geq Q(0)\wedge
\Big(\min_{t\in [0, \iy)}H(0,t)\Big) (=\del_1$ say) where we use the non-decreasing property of $\psi$. Hence $\p(\inf_{[0, T]}H(F^N(t), t)\geq\del_1/2)
\to 1$ as $N\to\iy$.

\skp

We note that $\chi$, defined in \eqref{inv}, may have singularity at $0$ if $y_{\max}\vee y_{0\max}=\iy$. Thus
$\chi$ is not a suitable function to be used in the continuous mapping theorem. Therefore we define
\[\tilde\chi(x, t)=\left\{
\begin{array}{ll}
\chi(x, t)\ & \text{for}\ x\geq\del_1/4, \ t\geq 0,
\\
\chi(\del_1/4, t) \ & \text{otherwise.}
\end{array}
\right.
\]
By definition, $\tilde\chi(\phi(\cdot), \cdot)=\chi(\phi(\cdot), \cdot)$.
Now we prove that $\tilde\chi$ is a continuous function.
Let $(x_n, t_n)\to(x, t)\in \R_+\times[0, T]$
and $\tilde\chi(x_n, t_n)=y_n$. We  show that $y_n\to y=\tilde\chi(x, t)$
as $n\to\iy$. $\chi$ being monotonically decreasing in $x$, we have
$y_n\in [0, \chi(\del_1/4, t_n)]$. Again $\sup_{[0, T]}\tilde\chi(\del_1/4, t)<\iy$.
Otherwise, if $z_n=\tilde\chi(\del_1/4, s_n)\to\iy$ for some sequence $\{s_n\}\subset [0, T]$, we have $0<\del_1/4=H(z_n,s_n)\to 0$
which is a contradiction. This implied that $\{y_n\}$ is bounded and hence have a subsequence, denoted again by $\{y_n\}$, converging to $y_0$.
We need to show that $y_0=y$.

\skp

\noi{Case 1:}\ $x\geq H(0, t)$ and so $y=0$.
 For the subsequence satisfying $x_n\geq H(0, t_n)$ , we have $y_n= 0$
 and the convergence is trivial.
 Hence we consider the subsequence
for which $x_n\in (\del_1/4, H(0, t_n))$. Then $H(y_n, t_n)=x_n\to H(0, t)$. By the convergence of $\{y_n\}$, $H(y_0, t)=H(0, t)$
and $y_0\in [0, \chi(\del_1/4, t)]$.
 Hence by strict monotonicity of $H(\cdot, t)$, $y_0=0$. Similar argument holds for the case
when $x\leq\del_1/4$.

\skp
\noi{Case 2:}\ $x\in(\del_1/4, H(0, t))$ and $x_n\in (\del_1/4, H(0, t_n))$. Hence
$x_n=H(y_n, t_n)\to x=H(y, t)$. Using continuity of $H$, we obtain $H(y_0, t)=H(y, t)$ and $y\in [0, \chi(\del_1/4, t)]$. This
implies $y=y_0$.

Define a map $\Psi: D_{\R_+}([0, \iy))\to D_{\R_+}([0, \iy))$
 by $\Psi(\zeta)(t)=\tilde\chi(\zeta(t), t)$. It is easy to see that $\Psi$ is continuous.
Hence applying the continuous mapping theorem,
 we have $\Psi(H(F^N(\cdot), \cdot))\Ra F$ as $N\to\iy$ and $F(t)=\chi(\phi(t), t)$. To complete the proof it is enough
 to show that $\sup_{[0, T]}|\Psi(H(F^N(t), t))-F^N(t)|\to 0$ in probability as $N\to\iy$. But this is obvious as on $\{\inf_{[0, T]}H(F^N(t), t)\geq\del_1/2\}$ and $\{\inf_{[0, T]}\frac{1}{N}Q^N(t)>0\}$, we have
 \[
 \Psi(H(F^N(t), t))=\tilde\chi(H(F^N(t), t), t)=F^N(t).
 \]

For general $y_{0\min}$, we note that by Assumption \ref{initial-meas}, $y_{0,\min}\leq y^*$. Hence for any $\kap>0$, $\tilde\chi(\cdot, \cdot)$
is continuous on $[\kap, T]$ for any $\kap>0$. Therefore $F^N\Ra F$ on $[\kap, T]$. \hfill $\Box$

\smallskip
Now we give the proof of our main result Theorem \ref{main-th}(b).

\medskip
\noi{\bf Proof of Theorem \ref{main-th}(b):}\ For $(x, t)\in\R_+\times\R_+$, define a measure $\calG(x,t)$ by
\begin{equation}\label{230}
\calG(x,t)(B)=\calQ_0(B\cap[x,\iy)+t)+\la\int_0^t\lan\Ir_{\{B\cap[x,\iy)+t-s\}}, \nu\ran ds,
\end{equation}
for $B\in\calB(\R_+)$. From Lemma \ref{lem-5}
and Corollary \ref{cor-1} we get that
\[
\sup_{a\in\R_+}\sup_{[0,T]}\Bigl|\frac{1}{N}\calQ^N_t(a, \iy)-\calG(F^N(t), t)(a,\iy)\Bigr|\to 0
\]
in probability. By our assumption on $\calQ_0$ and $\nu$, $\calG:\R_+\times\R_+\to\calM$ is a continuous map. To see this first we observe that for any $x\in\R_+$ we have
$$\lim_{\epsilon\to 0}\int_0^T\nu([x-\epsilon+s, x+\epsilon+s])ds=\int_0^T\nu(\{x+s\})ds=0,$$
where we used the dominated convergence theorem and then fact that $\nu$ has only countably many atoms.
Therefore if $(x_n, t_n)\to(x, t)$
then for any $f\in C_b(\R_+)$
\begin{align*}
\lan f, \calG(x_n, t_n)\ran&= \int_{[x_n+t_n, \iy)}f(y-t_n)\calQ_0(dy)+\int_0^{t_n}\int_{[x_n+t_n-s,\iy)}f(y-t_n+s)\nu(dy)ds
\\
&=\int_{[x_n+t_n, \iy)}f(y-t_n)\calQ_0(dy)+\int_0^{t_n}\int_{[x_n+s,\iy)}f(y-s)\nu(dy)ds
\\
&\rightarrow \int_{[x+t, \iy)}f(y-t)\calQ_0(dy)+\int_0^{t}\int_{[x+s,\iy)}f(y-s)\nu(dy)ds
\\
&=\int_{[x+t, \iy)}f(y-t)\calQ_0(dy)+\int_0^{t}\int_{[x+t-s,\iy)}f(y-t+s)\nu(dy)ds =\lan f, \calG(x, t)\ran
\end{align*}
as $n\to\iy$.
Hence
$(F^N(\cdot), \cdot)\Ra (F(\cdot), \cdot)$ implies that $\calG(F^N(\cdot), \cdot)\Ra \calG(F(\cdot), \cdot)$ as $N\to\iy$.
$F(\cdot)$ being continuous we in fact have, $\sup_{a\in\R_+}\sup_{[0,T]}|\calG(F(t),t)(a, \iy)-\calG(F^N(t), t)(a,\iy)|\to 0$
in probability as $N\to\iy$. Hence
$$ \sup_{a\in\R_+}\sup_{[0,T]}\Bigl|\frac{1}{N}\calQ^N_t(a, \iy)-\calG(F(t), t)(a,\iy)\Bigr|\to 0\quad \text{in probability},$$
as $N\to\iy$. Since $\calG(F(\cdot), \cdot)$ is continuous and non-atomic, it is easy to see that
$$ \sup_{[0,T]}\rho\Bigl(\frac{1}{N}\calQ^N_t,\calG(F(t), t)\Bigr)\to 0\quad \text{in probability},$$
as $N\to\iy$. Hence the proof follows.\hfill $\Box$

\subsection{Proof for the case $\la<\mu$}
In this section we give the proof for the sub-critical case.

\noi{\bf Proof of Theorem \ref{main3-th}:}\ To prove the theorem it is enough to show that for any $T, \eps>0$,
\begin{equation}\label{450}
\lim_{N\to\iy}\p\Bigl(\sup_{[0, T]}|\bar Q^N(s)-\bar\phi(s)|>\eps\Bigr)=0,
\end{equation}
where $\bar Q^N=\frac{1}{N}Q^N$. Recall time $\bar\tau^N$ from Lemma
\ref{lem-6},
$$\bar\tau^N=\inf\{t\geq 0\ :\ \bar Q^N(t)=0 \}. $$
Since the system observes a busy period in $[0, \tau^N)$, for
any $\eps_1, \del$ we have
\begin{equation}\label{451}
\lim_{N\to\iy}\p\Bigl(\sup_{[0, \bar T\wedge\tau^N]}\frac{1}{N}R^{N,\del}_1(s)\geq\eps_1\Bigr)=0.
\end{equation}
Following the same calculations as in Theorem \ref{main2-th}, we have
on $\Om_A(N)\cap\Om_S(N)\cap\Om_H(N)$,
\begin{align*}
&\sup_{[0, \bar T\wedge\tau^N]}|\bar Q^N(s)-\phi(s)|\leq\\
&
2(\eps_A(N)+\eps_S(N)+\eps_H(N))+4\del
 |\bar Q^N(0)-Q(0)|+\sup_{[0, T\wedge\tau^N]}\frac{1}{N}R^{N,\del}_1(s).
\end{align*}
Recall that $\bar T$ is the time when $\phi$ hits $0$ for the first time.
Without loss of generality we can assume that $\bar T\leq T$. Otherwise
the proof become obvious using the above inequality as $\p(\tau^N>T)\to 1$. Thus we assume $\bar T\leq T$.
Choose $\bar\del>0$ such that $\phi(\bar T-\bar\del)=\eps/4\wedge Q(0)$
and $\bar\phi(s)\leq\eps/4$ for $s\geq\bar T-\bar\del$ .
Now using the above inequality we can find events $\Om_{Q,\phi}(N)$ and a sequence $\{\eps_{Q, \phi}(N)\}, \eps_{Q, \phi}(N)\to 0,$ such that
$\p(\tilde\Om_{Q, \phi}(N))\geq 1-\eps_{Q, \phi}(N),$  and on $\Om_{Q,\phi}(N)$
\begin{equation}\label{455}
\sup_{[0, \bar T\wedge\tau^N]}|\bar Q^N(s)-\phi(s)|\leq \eps_{Q, \phi}(N),
\end{equation}
 This implies that on $\Om_{Q,\phi}(N)$,
${\tau^N>\bar T-\bar\del}$ for $N$ large. Again on $\Om_{Q,\phi}(N)$,
for large $N$ we have $\bar Q^N(\bar T-\bar\del)\leq \frac{3}{8}\eps$.
Hence for $s, t\in[\bar T-\bar\del, T]$, if $(s, t]$ is a busy period, then
\begin{align}\label{454}
\bar Q^N(t) &\leq \bar Q^N(s)+\bar E^N(s, t]
-\bar D^N(s, t]-(\bar R^N(t)-\bar R^N(s))\notag
\\
&\leq \bar Q^N(s)+2(\eps_A(N)+\eps_S(N))+\la(t-s)-\mu\int_s^t\Ir_{\bar X^N(u)>0}du,
\end{align}
on $\Om_A(N)\cap\Om_S(N)$, where in the last inequality we use the monotone
property of $R^N(\cdot)$. Thus if we choose $N$ large so that
$2(\eps_A(N)+\eps_S(N))<\eps/8$, then
$\sup_{[\bar T-\bar\del, T]}\bar Q^N\leq\eps/2$ on $\Om_A(N)\cap\Om_S(N)\cap\Om_{Q,\phi}(N)$. To see this observe that if
there is an interval $[s, t]\subset[\bar T-\bar\del, T]$ so that
$\bar Q^N(s)\leq 3\eps/8, \bar Q^N(t)> \eps/2$ and $\bar Q^N>0$
on $(s, t]$, we shall have a contradiction from \eqref{454} as $\la<\mu$.

Hence on $\Om_A(N)\cap\Om_S(N)\cap\Om_{Q,\phi}(N)$,
\begin{align*}
\sup_{[0, T]}|\bar Q^N(s)-\bar\phi(s)| &\leq \sup_{[0, \bar T-\bar\del]}|\bar Q^N(s)-\bar\phi(s)|+ \sup_{[\bar T-\bar\del, T]}|\bar Q^N(s)-\bar\phi(s)|
\\
&\leq \eps/8+\eps/2+\eps/4<\eps,
\end{align*}
where in the last inequality we used \eqref{455}. Hence the proof.\hfill $\Box$

\noindent{\bf Acknowledgement.}
We would like to thank the referee for many helpful comments.

\footnotesize

\bibliographystyle{is-abbrv}

\bibliography{refs}

\iffalse

\fi

\end{document}